\documentclass[11pt]{article}

\usepackage{geometry}
\date{}

\usepackage{amsthm}
\usepackage{graphics} 
\usepackage{epsfig} 
\usepackage{times}
\usepackage{amsmath} 
\usepackage{amssymb}  
\usepackage{color}
\newtheorem{assumption}{Assumption}
\newtheorem{lemma}{Lemma}

\newtheorem{remark}{Remark}
\newtheorem{definition}{Definition}

\usepackage{bbm}
 \usepackage{url}
\usepackage{amssymb} 
\usepackage{amsmath}
\DeclareMathOperator*{\argmin}{argmin}

\usepackage[utf8]{inputenc} % allow utf-8 input
\usepackage[T1]{fontenc}    % use 8-bit T1 fonts
\usepackage{hyperref}       % hyperlinks
\usepackage{url}            % simple URL typesetting
\usepackage{booktabs}       % professional-quality tables
\usepackage{amsfonts}       % blackboard math symbols
\usepackage{nicefrac}       % compact symbols for 1/2, etc.
\usepackage{microtype}      % microtypography
\usepackage{xcolor}         % colors
\usepackage[font=small,labelfont=bf]{caption}               % captions
\usepackage{wrapfig}
\usepackage{float}
\usepackage{algorithm}
\usepackage{algpseudocode}
\usepackage{thmtools, thm-restate}

\usepackage{caption}
\usepackage{subcaption}
\usepackage{enumitem}

\usepackage{ifthen}
\newif\ifshortver
\usepackage{caption}

\hypersetup{
    bookmarks=true,         % show bookmarks bar?
    unicode=false,          % non-Latin characters in Acrobat’s bookmarks
    pdftoolbar=true,        % show Acrobat’s toolbar?
    pdfmenubar=true,        % show Acrobat’s menu?
    pdffitwindow=false,     % window fit to page when opened
    pdfstartview={FitH},    % fits the width of the page to the window
    pdfnewwindow=true,      % links in new PDF window
    colorlinks=true,       % false: boxed links; true: colored links
    linkcolor=blue,          % color of internal links (change box color with linkbordercolor)
    citecolor=blue,        % color of links to bibliography
    filecolor=magenta,      % color of file links
    urlcolor=blue,           % color of external links
    breaklinks=true
}

\newcommand{\mc}{\mathbb}

\usepackage{tikz}

\title{Coreset-Based Task Selection for Sample-Efficient\\ Meta-Reinforcement Learning}

\usepackage[round]{natbib}

\allowdisplaybreaks
\usepackage[toc,page,header]{appendix}
\usepackage{minitoc}

\usepackage{amssymb} 
\usepackage{amsmath}

\usepackage{graphicx}
\usepackage{float}
\usepackage{algorithm}
\usepackage{algpseudocode}
\usepackage{listings,lstautogobble}
\usepackage{tikz}
\usepackage{courier}
\usetikzlibrary{decorations.markings}
\usepackage{graphicx}
\usepackage{color}
\usepackage{hyperref}
\usepackage{cleveref}
\usepackage{caption}
\usepackage{bbm, dsfont}
\usepackage{wrapfig}

\newcounter{myvar}
\newcommand\blfootnote[1]{%
  \begingroup
  \renewcommand\thefootnote{}\footnote{#1}%
  \addtocounter{footnote}{-1}%
  \endgroup
}

\newtheorem{cor}{Corollary}

\newtheorem{thrm}{Theorem}

\allowdisplaybreaks

\author{Donglin  Zhan, Leonardo F. Toso, and James Anderson$^*$}

\begin{document}

\doparttoc % Tell to minitoc to generate a toc for the parts
\faketableofcontents % Run a fake tableofcontents command for the partocs

%\part{} % Start the document part
%\parttoc % Insert the document TOC

\maketitle

\begin{abstract}%
We study task selection to enhance sample efficiency in model-agnostic meta-reinforcement learning (MAML-RL). Traditional meta-RL typically assumes that all available tasks are equally important, which can lead to task redundancy when they share significant similarities. To address this, we propose a coreset-based task selection approach that selects a weighted subset of tasks based on how diverse they are in  gradient space, prioritizing the most informative and diverse tasks. Such task selection reduces the number of samples needed to find an $\epsilon$-close stationary solution by a factor of $\mathcal{O}(1/\epsilon)$. Consequently, it guarantees a faster adaptation to unseen tasks while focusing training on the most relevant tasks. As a case study, we incorporate task selection to MAML-LQR \citep{toso2024meta}, and prove a sample complexity reduction proportional to $\mathcal{O}(\log(1/\epsilon))$ when the task-specific cost also satisfy gradient dominance. Our theoretical guarantees underscore task selection as a key component for scalable and sample-efficient meta-RL. We numerically validate this trend across multiple RL benchmark problems, illustrating the benefits of task selection beyond the LQR baseline.
\end{abstract}

\section{Introduction}
\blfootnote{$^*$All authors are with the Department of Electrical Engineering, Columbia University in
the City of New York. Email: \texttt{\{dz2478, lt2879, james.anderson\}@columbia.edu}.}

\begin{wrapfigure}{r}{0.6\textwidth}
    \vspace{-0.28in}
    \centering
    \includegraphics[width=0.55\textwidth]{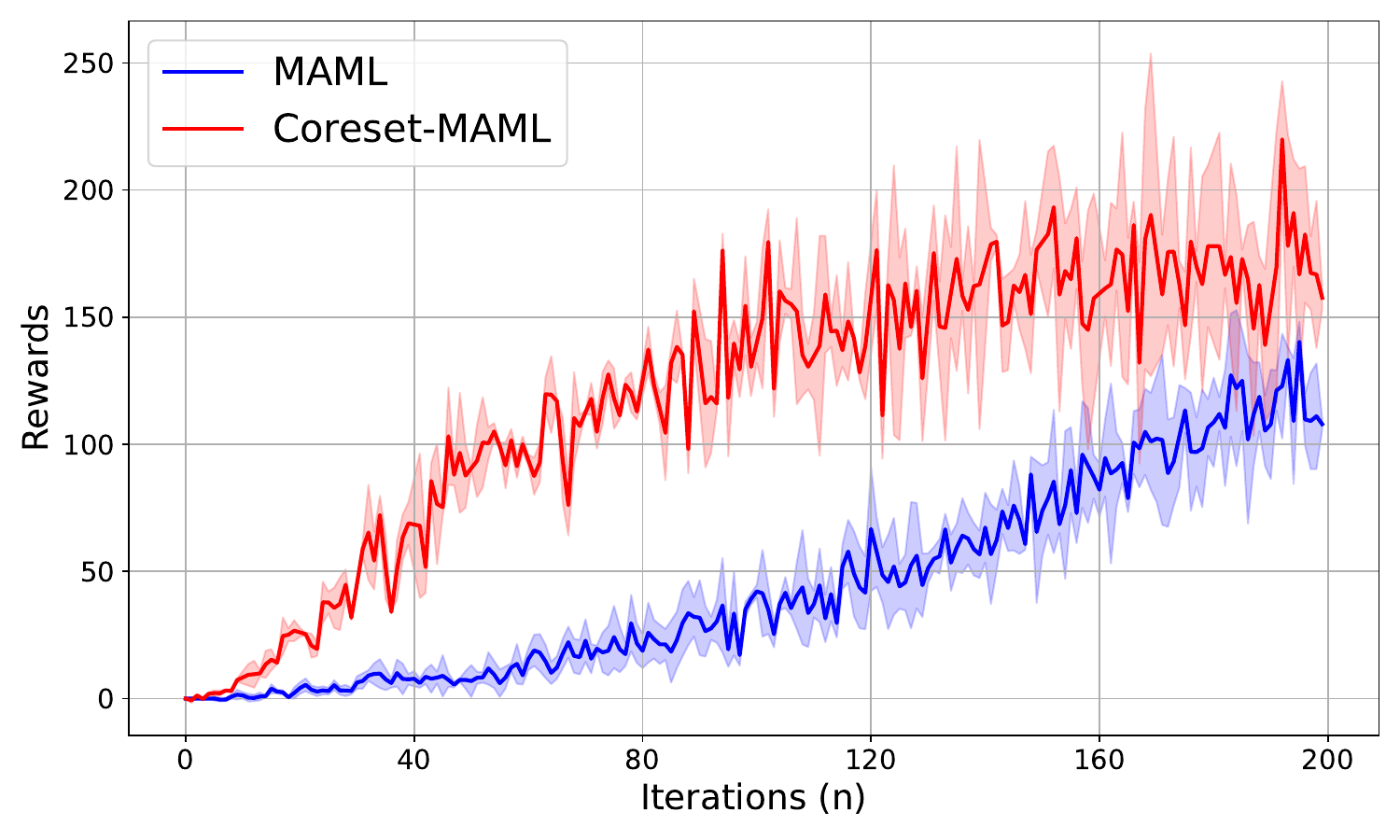}
    \vspace{-0.4cm}
    \caption{Comparison of coreset MAML-RL (this work) and  MAML-RL on the walker2D Mujoco environment (see Section \ref{sec:numerics} for details).}
    \label{fig:hopper-intro}
    \vspace{-0.1in}
\end{wrapfigure}

\noindent Meta-reinforcement learning (meta-RL) has emerged as a powerful framework for learning policies that can quickly adapt to unseen environments \citep{wang2016learning,finn2017model}. In particular, the model-agnostic meta-reinforcement learning (MAML-RL) algorithm has demonstrated success in enabling agents to learn a shared policy initialization that is only a few policy gradient steps away from optimality for any  seen \emph{and} unseen task \citep{duan2016rl,nagabandi2018learning}. Such quick adaptation is crucial, for example, in robotics \citep{song2020rapidly}, where agents often need to operate in dynamic environments  and accomplish a variety of  goals.

MAML-RL and meta-reinforcement learning more generally, typically assumes that all training tasks are equally important. This assumption may lead to task redundancy and excessive sampling costs as it is likely not worth sampling from multiple similar tasks; instead collecting data from a single representative task would suffice.

``Task selection'' can be thought of a pre-processing step in the meta-learning pipeline. It seeks to identify a representative subset of tasks that captures the diversity across all  training tasks, and then uses this smaller ``coreset'' for training. In particular, ``coreset learning'' has been proposed for data-efficient training of machine learning models \citep{mirzasoleiman2020coresets1, pooladzandi2022adaptive,yang2023towards}. Related work has also employed coreset selection to select clients in federated learning \citep{balakrishnan2022diverse} and continual learning \citep{tiwari2022gcr,wang2022learning}. For meta-learning, and in the context of classification, \cite{zhan2024data} propose a data-efficient and robust task selection algorithm (\texttt{DERTS}) that outperforms existing sampling-based techniques. In essence, \texttt{DERTS} frames coreset learning as submodular optimization, where the goal is to select a subset of tasks that minimizes the maximum normed difference across task-specific gradients. 

In the work above, it is   assumed that task-specific gradients can be directly computed. In meta-RL, such an assumption may be restrictive as the meta-gradient depends on unknown task and trajectory distributions, where the later is also conditioned on the current policy. As such, it is  challenging to compute the gradient via automatic differentiation \citep{rothfuss2018promp}. To circumvent this, one must resort to gradient approximation, that by itself introduces extra difficulties to the analysis of meta-RL with task selection.  Specifically,  errors arising from  gradient estimation and the meta-training on the coreset need to be carefully accounted for.

\vspace{0.1cm}
\noindent \textbf{Contributions:} Towards addressing these points, we propose a coreset-based task selection algorithm (inspired by \texttt{DERTS}) for meta-RL. The main contributions of our approach are:

\vspace{0.1cm}
\noindent $\bullet$ \textbf{Algorithmic:} This is the first work to propose a derivative-free coreset-based task selection approach for MAML-RL (Algorithm \ref{alg:metatrain}), which also comes with strong convergence guarantees  (Section \ref{sec:theory}). We derive an ergodic convergence rate for non-concave task-specific reward functions (Theorem \ref{theorem:ergodic convergence}) and prove that Algorithm \ref{alg:metatrain} finds an $\epsilon$-close stationary solution after $N = \mathcal{O}(1/\epsilon)$ iterations when the task selection bias is made sufficiently small. We also incorporate task selection to meta-learning for control via the \texttt{MAML-LQR} algorithm~\citep{toso2024meta} and show that it learns a provably fast-to-adapt LQR controller (Theorem \ref{theorem: global convergence guarantees}) within $N = \mathcal{O}(\log(1/\epsilon))$ iterations while reducing task redundancy.
\vspace{0.2cm}

\noindent $\bullet$ \textbf{Sample Complexity:} We demonstrate that  selecting a weighted subset of the most informative tasks reduces the sample complexity for achieving local convergence by a factor of $\textcolor{blue}{\mathcal{O}(1/\epsilon)}$ (c.f. Section~\ref{sec:numerics}, Figure \ref{fig:hopper-intro}, and Corollary \ref{corollary: sample complexity local}). In particular, this reduction is guaranteed when the set of training tasks is sufficiently large and tasks therein are sufficiently similar. Moreover, Algorithm \ref{alg:metatrain} offers a sample complexity reduction proportional to \textcolor{blue}{$\mathcal{O}(\log(1/\epsilon))$}  in the \texttt{MAML-LQR} setting (Corollary \ref{corollary: sample complexity global}).

\vspace{0.2cm}
\noindent \textbf{Related Work:} Meta-reinforcement learning has been extensively studied across several applications, including robot manipulation \citep{yu2017preparing}, locomotion \citep{song2020rapidly}, and building energy control \citep{luna2020information}. Most relevant to our work is \cite{song2019maml, song2020rapidly}, which is derivative-free but treats all tasks equally, leading to task redundancy. Task weighting is addressed in \cite{shin2023task} and \cite{zhan2024data} by selecting representative task subsets. However, \cite{zhan2024data} focuses on classification tasks and simplifies gradient approximation by using the model pre-activation outputs, while \cite{shin2023task} employs an information-theoretic metric for task selection and does not consider gradient-based training.  Finally, in the context of control, the linear quadratic regulator (LQR) problem has become a key baseline for policy optimization in reinforcement learning.  In particular, \cite{molybog2021does,musavi2023convergence,toso2024meta,aravind2024moreau,pan2024model} study the meta-LQR problem and provide guarantees for provably learning fast-to-adapt LQR controllers. Building on this, our work also integrates task selection into the MAML-LQR setting, demonstrating its effectiveness in reducing the sample complexity. A broader overview of related work is included in Section~\ref{sec:related} of the appendix.

\section{Preliminaries}

We now introduce the MAML-RL problem \citep{finn2017model} and formalize our coreset-based task selection that selects a weighted subset of tasks based on their diversity in  gradient space.

\subsection{Model Agnostic Meta-Reinforcement Learning}

Let $T_j$ be a reinforcement learning task drawn from a task distribution $\mathcal{P}(\mathcal{T})$ over a set of tasks $\mathcal{T}$. Let $\mathcal M$ denote a set of $M$ tasks drawn from $\mathcal{P}(\mathcal{T})$, i.e., $\mathcal{M}:= \left\{T_{j} \mid j = 1,2,\ldots, M\right\}$.  The objective of MAML-RL is to learn a meta-policy $\pi_{\theta^{\star}}$ with meta-parameter $\theta^\star \in \Theta$, trained on $\mathcal M$,  such that within a few policy gradient steps, $\pi_{\theta^{\star}}$ can be adapted to an unseen task-specific policy.  We assume that all tasks share the same state and action spaces $\mathcal{S}$, and $\mathcal{A}$.  In addition, each task $T_{j}$ is associated with a reward function $R_j$ and a transition distribution $q_{j}(s_{t+1}|s_t,a_t)$ at time step $t$. Also, let $J(\theta)$ and $J_{j}(\theta)$ be the MAML and task-specific reward functions, respectively. The one-shot MAML-RL problem is written as follows:
\begin{align}\label{eq:MAML_RL}
    \theta^\star := \text{argmax}_{\theta \in \Theta}~ J(\theta) := \mc{E}_{T_j \sim \mathcal{P}(\mathcal{T})} J_j(\theta + \eta_{\text{inn}} \nabla J_j(\theta)),
\end{align}
for some positive inner step-size $\eta_{\text{inn}}$. Moreover, $J_{j}({\theta}) := \mc{E}_{{\tau} \sim \mathcal{P}_{T_j} ({\tau} | {\theta})} \mathcal{R}_j({\tau})$ is the task-specific reward incurred by $\pi_{\theta}$, where $\mathcal{P}_{T_j}(\tau | \theta)$ is the distribution of trajectories $\tau$ conditioned on the policy $\pi_\theta$. We also define $\bar{\theta} :=  \theta + \eta_{\text{inn}} \nabla_\theta J_j(\theta)$. Hence, the gradient-based MAML-RL update follows: $\theta \leftarrow \theta + \eta_{\text{out}}\nabla J(\theta)$ with $\eta_{\text{out}}$ denoting the outer step-size, and the MAML gradient given by
\begin{align*}%\label{eq:MAML_RL_Grad}
\nabla J(\theta) = \mathbb{E}_{T_j\sim\mathcal{P}(\mathcal{T})}\left[\mathbb{E}_{{\tau}\sim\mathcal{P}_{T_j}({\tau}|\bar{\theta})}(\nabla_{\bar{\theta}} \log \mathcal{P}_{T_j}({\tau}|\bar{\theta})R_j({\tau})\nabla \bar{\theta})\right],  \text{ where}%\tag{2}
\end{align*}
$$
\nabla_{\bar{\theta}} := I + \eta_{\text{out}}\left( \int \mathcal{P}_{T_j}(\tau|\theta)\nabla^2\log \pi_\theta(\tau)R_j(\tau)+\mathcal{P}_{T_j}(\tau|\theta)\nabla_\theta\log \pi_\theta(\tau)\nabla_\theta\log \pi_\theta(\tau)^\top R_j(\tau)d\tau\right). %\tag{3}
$$
The direct computation of $J_j(\theta)$ and $\nabla J(\theta)$ may be intractable due to the expectation over unknown trajectory and task distributions. A standard approach to approximate $J_j(\theta)$ is to sample multiple trajectories $\tau_l$ by playing  $\pi_{\theta}$ and then computing the empirical task-specific reward function $\hat{J}_j(\theta)$ \citep{liu2019taming,rothfuss2018promp,song2019maml}. On the other hand, following \cite{song2019maml,song2020rapidly}, we propose a derivative-free method for estimating the task-specific and MAML gradients through querying/estimating task-specific rewards.

\subsection{Task Selection}
Motivated by \cite{yang2023towards, zhan2024data}, 
we argue that not all the tasks in the task pool $\mathcal{M}$ are equally important for  meta-training. Since multiple tasks may share similarities, that may lead to task redundancy and sample inefficiency, as it requires collecting data from multiple similar tasks when, in principle, collecting data from a single  task that is representative of all the similar tasks
should  suffice for training. Our goal is to select a subset of tasks, $\mathcal S$ (the coreset), from the task pool $\mathcal M$ that best represents the diversity of the  tasks in  $\mathcal{M}$, such that the performance of the model trained on a weighted subset of $\mathcal S$  is  sufficiently close to that of a model trained on the full  task pool $\mathcal M$. In particular, we will prove that by carefully controlling the task selection bias one may achieve a substantial sample complexity reduction in the meta-training.

The main steps of our meta training algorithm are: $(i)$ coreset selection---selecting the coreset $\mathcal S \subseteq \mathcal{M}$, $(ii)$ weight allocation for  each task in $\mathcal S$ such that it captures the relative  importance of the task, and $(iii)$  meta-training on the coreset $\mathcal S$ and its corresponding weighting. In the following subsections we  introduce the concept of gradient approximation  over the task pools and establish the selection criterion for assigning tasks in $\mathcal M$ to $\mathcal S$.

\subsubsection{Full Gradient Approximation over the Task Pool} \label{sub1}

Our aim is to select the coreset $\mathcal{S} \subseteq \mathcal{M}$ with $\mathcal{S} := \left\{T_{i} \mid i = \alpha_1,\alpha_2,\ldots,\alpha_L\right\}$, where $\alpha_i \in [M]$ and $L \leq M$, with corresponding weights $\left\{\gamma_{i} ~|~ i = 1,2,\ldots,L\right\}$, such that the gradient for training on $\mathcal S$ with corresponding weights $\gamma_{i}$ approximates the meta-gradient on $\mathcal{M}$. 

To better understand the coreset-based task selection we let $\Gamma: \mathcal{M} \rightarrow \mathcal{S}$ be a mapping from the task pool $\mathcal{M}$ to the coreset $\mathcal{S}$, i.e., that maps a task $T_{j}$ from $\mathcal{M}$ to a task $T_{i}$ in $\mathcal{S}$. For simplicity, we denote $\Gamma(T_{j}) = T_{i}$ as $\Gamma(j) = i$. In addition, let $\mathcal{S}_c$ denote the complement of $\mathcal{S}$ in $\mathcal{M}$.  Following \cite{mirzasoleiman2020coresets1}, we define the weight $\gamma_i$ of the selected task $T_{i} \in \mathcal{S}$ as
$\gamma_i := \sum_{j \in \mathcal{M}}\mathbbm{1}_{\{\Gamma(j)=i\}}$, where $\mathbbm{1}_\mathcal C$ is the indicator function over some set $\mathcal C$. Also, the summation over $\mathcal{M}$ results the number of  tasks in $\mathcal{M}$ that are assigned to task $T_{i}$ in $\mathcal{S}$.

Then, by using the definition of the mapping function $\Gamma$ the gradient approximation error due to training over $\mathcal{S}$ instead of $\mathcal{M}$ is given by 
\begin{equation}
\begin{split}
\label{approx}
\hspace{-0.1cm}\left\|\sum_{j \in \mathcal{M}} \nabla J
_{j}(\theta) - \sum_{j \in \mathcal{M}} \nabla J
_{\Gamma(j)}(\theta)\right\| 
\hspace{-0.1cm}= \hspace{-0.1cm}\left\|\sum_{j \in \mathcal{M}} \nabla J
_{j}(\theta) - \sum_{i \in \mathcal{S}} \gamma_{i} \nabla J_{i}(\theta)\right\| \hspace{-0.1cm} \leq \hspace{-0.1cm} 
\sum_{j \in \mathcal{M}} \min _{i \in \mathcal{S}}\left\| \nabla J
_{j}(\theta) - \nabla J_{i}(\theta)\right\|,
\end{split}
\end{equation}
where our objective is to control and make this error as small as possible (by selecting $\mathcal S$ and the weights). We emphasize that we \emph{do not} have access to $\mathcal{S}$ and subsequently $\Gamma$. Namely, we cannot directly evaluate that error, and optimizing over the subset of tasks $\mathcal S$ is NP-hard. Instead, we proceed by minimizing the  RHS of~\eqref{approx}. Namely, by assuming the elements in $\mathcal{S}$ are fixed we assign each task in $\mathcal{M}$ to its closest element in $\mathcal{S}$, in the gradient space, through the mapping $\Gamma$. To do so, the weights $\gamma_{i}$, for all tasks $T_i \in \mathcal{S}$, corresponding to the mapping $\Gamma$ can be allocated as $
\gamma_i=\sum_{j \in \mathcal{M}} \mathbbm{1}_{\left\{j=\operatorname{argmin}_{T_{i} \in \mathcal{S}}\left\| \nabla J
_{j}(\theta) - \nabla J_{i}(\theta)\right\|\right\}}.$

However, as previously discussed, directly computing $\nabla J_i(\theta)$ may not be tractable for most RL tasks. This is due to the fact that simulating trajectories and performing backpropagation through deep-RL models incur a large computational cost. This motivates the use of derivative-free methods to approximate the gradient. In particular, we propose a derivative-free task selection approach based on a zeroth-order gradient estimation scheme. 

For this, we consider a two-point estimation since it has a lower estimation variance compared to its one-point counterpart \citep{malik2019derivative}.  Zeroth-order estimation is a Gaussian smoothing approach \citep{nesterov2017random} based on Stein's identity \citep{stein1972bound} that relates gradient to reward queries. We refer the reader to \cite{flaxman2004online,spall2005introduction} for further details on zeroth-order gradient estimation. The two-point zeroth-order estimation of $\nabla J_{i}(\theta)$ is
\begin{align}\label{eq:ZO2P}
    \texttt{ZO2P}(\theta,n_s,r) :=  \widehat{\nabla} J_i(\theta) := \frac{d}{2n_sr^2}\sum_{l = 1}^{n_s}\left(J_{i}(\theta + v_l) - J_i(\theta - v_l)\right)v_l,
\end{align}
with $r > 0$ denoting the smoothing radius, $v_l \in \mathbb{R}^{d_1\times d_2}$ randomly drawn from a uniform distribution over the Euclidean sphere of radius 
$r$, $\mathbb{S}^{d-1}_r$, namely, $v_l \sim \mathbb{S}^{d-1}_r$, $n_s$ being the number of samples, and $d := d_1d_2$. Finally, the estimation of the meta-gradient over $\mathcal{M}$ and $\mathcal{S}$ is given by
\begin{align*}
    \nabla_{\mathcal{M}} J(\theta) := \frac{1}{M}\sum_{j \in \mathcal{M}} g_j(\theta), \;\ \nabla_{\mathcal{S}} J(\theta) := \frac{1}{M}\sum_{i \in \mathcal{S}}\gamma_i g_i(\theta),
\end{align*}
where $g_i(\theta) := \frac{d}{2r^2n_s} \sum_{l = 1}^{n_s} \left(J_{i}(\theta + u_l + \eta_{\text{inn}}\widehat{\nabla} J_i(\theta)) - J_i(\theta - u_l + \eta_{\text{inn}}\widehat{\nabla} J_i(\theta))\right)u_l$.

\subsubsection{Coreset-based Task Selection} \label{sub2}

We note that minimizing the RHS of~\eqref{approx} is equivalent to maximizing the ``facility location'' function, which is a well-known submodular function \citep{cornuejols1977uncapacitated}. 

\begin{definition}[Submodularity~\citep{nemhauser1978analysis}]\label{submod}
A set function $F: 2^{V} \rightarrow \mathbb{R}^{+}$ is submodular if $F(e \mid S) := F(S \cup\{e\})- F(S) \geq F(T \cup\{e\})-F(T)$, for any $S \subseteq T \subseteq V$ and $e \in V \backslash T$. $F$ is monotone if $F(e \mid S) \geq 0$ for any $e \in V \backslash \bar{S}$ and $S \subseteq V$. 
\end{definition}

We  leverage Definition~\ref{submod} and \eqref{approx} to define a monotone submodular function $\mathcal{F}$ over $\mathcal{S}$ with respect to the zeroth-order approximated gradient $g_{i}(\theta)$. That is,
\begin{equation}\label{eq:submod}
\mathcal{F}(\mathcal{S}):= C -\sum_{j \in \mathcal{M}} \underset{i \in \mathcal{S}}{\min}~\left\| g_{j}(\theta) -  g_{i}(\theta)\right\|,
\end{equation}
where $C>0$ upper bounds $\mathcal{F}(\mathcal{S})$. Therefore, to formulate our task selection objective, we restrict the cardinality of $\mathcal{S}$ and make the number of tasks in $\mathcal{S}$ sufficiently small ($L \ll M$). This is introduced as the constraint $|\mathcal{S}| \leq L$ in the following submodular optimization. The coreset $\mathcal S$ is learned by solving
\begin{equation}\label{eq:submodular_opt}
    \mathcal{S}^{\star} =\underset{\mathcal{S} \subseteq \mathcal{M}}{\arg \max }~ \mathcal{F}(\mathcal{S}), \text { s.t. }|\mathcal{S}| \leq L. 
\end{equation}

It is well-known that \eqref{eq:submodular_opt} can be solved through a greedy-based approach with a $1-e^{-1}$ error bound on the corresponding approximate solution \citep{nemhauser1978analysis,wolsey1982analysis}. We then incorporate such coreset-based task selection in the MAML-RL training over the learned weighted subset $\mathcal{S}$ in Algorithm \ref{alg:metatrain}. To start with, we initialize $\mathcal{S}$ as the empty set in step 1, and for each greedy iteration, we select a task $T_i$ from $\mathcal{S}_c$ that maximizes the marginal utility $\mathcal{F}(T_i|\mathcal{S}) = \mathcal{F}(\mathcal{S} \cup T_i) - \mathcal{F}(\mathcal{S})$ in steps 7 and 8, and update $\mathcal{S}$ as $\mathcal{S} = \mathcal{S} \cup {\arg\max}_{T_{i} \in \mathcal{S}_c} \mathcal{F}(T_{i} \mid \mathcal{S})$. With the learned subset $\mathcal{S}$ in hand, the weights $\gamma_i$ for all tasks $T_i \in \mathcal{S}$ are allocated in step 10. The task selection is then followed by step 12 to 17 where MAML is applied on the coreset $\mathcal S$.  In the next section we present the theoretical guarantees of Algorithm \ref{alg:metatrain}.

\begin{algorithm}
\caption{\texttt{Coreset Selection for MAML-RL}}
\label{alg:metatrain}
\begin{algorithmic}[1]
    \State \textbf{Input:} initial meta-policy parameter $\theta_0$; step-sizes $\eta_{\text{inn}}, \eta_{\text{out}}$; number of samples $n_s$; smoothing radius $r$; number of iterations $N$; number of selected tasks $L$; task pool $\mathcal{M}$; coreset $\mathcal{S} = \emptyset$
    \State \textcolor{blue}{\textbf{\texttt{Coreset Selection:}}}
    \State \textbf{for all} tasks $T_j$ in $\mathcal{M}$ \textbf{do} 
    \State \quad  $\widehat{\nabla} J_j(\theta_0) \leftarrow \texttt{ZO2P}(\theta_0,n_s,r)$,~  $g_j(\theta_0) \leftarrow \texttt{ZO2P}(\theta_0 + \eta_{\text{inn}}\widehat{\nabla} J_j(\theta_0),n_s,r)$ \Comment{\textcolor{red}{\texttt{estimation}}}
    \State \textbf{end for}
    \State \textbf{while} $|\mathcal{S}| < L$ \textbf{do}
    \State \quad $T_{i} \in \text{argmax} _{T_i \in S^{C}}\mathcal{F}(T_{i} \mid \mathcal{S})$ 
    \State \quad $\mathcal{S}=\mathcal{S} \cup\{T_{i} \}$
    \State  \textbf{end while}
    \State $\gamma_i=\sum_{j \in \mathcal{M}} \mathbbm{1}_{\left\{j=\operatorname{argmin}_{T_{i} \in \mathcal{S}}\left\| g_j(\theta_0) - g_i(\theta_0) \right\|\right\}}$ \Comment{\textcolor{red}{\texttt{weight allocation}}}
    \State \textcolor{blue}{\textbf{\texttt{MAML-RL over $\mathcal{S}$:}}}
    \State \textbf{for all} iteration $n = \{0,1,\ldots, N-1\}$ \textbf{do}
        \State \quad \textbf{for all} tasks $T_i$ in $\mathcal{S}$ \textbf{do}
            \State \quad \quad $\widehat{\nabla} J_i(\theta_n) \leftarrow \texttt{ZO2P}(\theta_n,n_s,r)$, ~ $g_i \leftarrow \texttt{ZO2P}(\theta_n + \eta_{\text{inn}}\widehat{\nabla} J_i(\theta_n),n_s,r)$ \Comment{\textcolor{red}{\texttt{estimation}}}
        \State \quad \textbf{end for}
        \State \quad  $\nabla_{\mathcal{S}}J(\theta_n) = 
        \frac{1}{M}\sum_{i \in \mathcal{S}} \gamma_i g_i(\theta_n), \:\ \theta_{n+1} = \theta_{n} + \eta_{\text{out}}\nabla_{\mathcal{S}}J(\theta_n)$ \Comment{\textcolor{red}{\texttt{meta update}}}
    \State \textbf{end for}
\State \textbf{Output:} $\theta_N$
\end{algorithmic}
\end{algorithm}

\section{Theoretical Guarantees} \label{sec:theory} %Leo

We first introduce the ergodic convergence rate (i.e., local convergence analysis) of Algorithm~\ref{alg:metatrain} for the general case of non-concave task-specific reward function. We then extend our results to meta-learning for control, specifically by applying task selection to the \texttt{MAML-LQR} algorithm from \cite{toso2024meta}, and derive global convergence guarantees when the task-specific cost satisfies a gradient dominance property. In addition, we discuss the sample complexity reduction benefit of task selection for both MAML-RL and MAML-LQR problems.

\subsection{Ergodic Convergence Rate}

For the local convergence analysis, when $J_i(\theta)$ is generally non-concave, our goal is to characterize the ergodic convergence rate with respect to the MAML reward function \eqref{eq:MAML_RL}, namely, we aim to control $\frac{1}{N}\sum_{n = 0}^{N-1}\|\nabla J(\theta_n)\|^2$. We first assume that the task-specific reward function and its gradient are locally smooth and that the gradient is  uniformly upper bounded. In addition, we assume that the upper bound of $\mathcal{F}(\mathcal{S})$ can be made sufficiently small.

\begin{assumption} (Local smoothness) \label{assump: Lipschitz of the MAML task specific} The task-specific reward function $J_i(\theta)$ and its gradient $\nabla J_i(\theta)$ are smooth with constants $\beta$ and $\psi$, respectively, i.e., for for any $\theta, \theta^\prime \in \Theta$, we have
\begin{align}\label{eq: assump Lipschitz}
    &|J_i(\theta) - J_i(\theta^\prime)| \leq \beta J_i(\theta) \|\theta - \theta^\prime\|,\;\ \|\nabla J_i(\theta) - \nabla J_i(\theta^\prime)\| \leq \psi \|\theta - \theta^\prime\|.
\end{align} 
\end{assumption}
\begin{assumption} (Gradient uniform bound) \label{assump: uniform bound task specific grad} $\|\nabla J_i(\theta)\| \leq \phi$, for any $T_i \in \mathcal{P}(\mathcal{T})$ and $\theta \in \Theta$.
\end{assumption}
\vspace{-0.2cm}

Assumptions \ref{assump: Lipschitz of the MAML task specific} and \ref{assump: uniform bound task specific grad} are standard in the convergence analysis of training dynamics \citep{oymak2019overparameterized,liu2022loss}, as well as in the literature of stochastic gradient descent (SGD) for non-convex loss functions \citep{stich2018local,li2019convergence}. Later, for the MAML-LQR setting, such conditions are in fact properties of the task-specific LQR cost. 

\begin{assumption} \label{assumption: submodular function upper bound} The constant $C$ in~\eqref{eq:submod} is set sufficiently small, i.e., $C = \mathcal{O}(\epsilon)$, for some small $\epsilon$. 
\end{assumption}

It is worth noting that making $C$ sufficiently small is standard in coreset learning for data-efficient machine-learning \citep{mirzasoleiman2020coresets1, yang2023towards} as it guarantees that task selection estimation error remains sufficiently small. Next, we define the maximum normed difference between the gradient of two task-specific reward functions over the parameter space and present the local convergence guarantee of Algorithm \ref{alg:metatrain}.

\begin{definition}\label{definition: maximum normed difference} The maximum normed difference between two distinct task-specific gradients over the parameter space $\theta \in \Theta$ is $\xi_{i,j} := \max_{\theta \in \Theta}\|\nabla J_i(\theta) - \nabla J_j(\theta)\|$.
\end{definition}

\begin{thrm} (Stationary solution)\label{theorem:ergodic convergence} Suppose Assumptions \ref{assump: Lipschitz of the MAML task specific}, \ref{assump: uniform bound task specific grad} and \ref{assumption: submodular function upper bound} are satisfied. In addition, suppose the number of samples and smoothing radius are set according to 
\begin{align*}
n_s \geq C_{\text{approx},1}\left(\frac{\sigma^2}{\epsilon^2}+\frac{ b}{3\epsilon}\right) \log \left(\frac{d_1 + d_2}{\delta}\right),\quad r \leq \frac{\epsilon}{C_{\text{approx},1}\psi},
\end{align*}
with $\sigma^2 :=  \left(d\beta J_{\max}\right)^2 + \left(\epsilon + \phi\right)^2$ and $b := d\beta J_{\max} +   \epsilon +\phi$ for $J_{\max}:= \max_{T_j \in \mathcal{M}, \theta \in \Theta} J_j(\theta)$, $\delta, \epsilon \in (0,1)$, for a sufficiently large universal constant $C_{\text{approx},1}$. Lastly, suppose that the step-sizes scale as $\eta_{\text{inn}} = \mathcal{O}(\epsilon/d)$ and $\eta_{\text{out}} = \mathcal{O}(1)$, and $L = \mathcal{O}(1)$. Then, Algorithm \ref{alg:metatrain} satisfies
\begin{align*}
    \frac{1}{N}\sum_{n = 0}^{N-1}\|\nabla J(\theta_n)\|^2_2 &\leq   \mathcal{O}\left(\frac{\Delta_0}{\eta_{\text{out}}N}  +  \left(\frac{1}{M}\sum_{j \in \mathcal{M}} \min_{i \in \mathcal{S^\star}}\xi_{i,j}\right)^2\right),
\end{align*}
with probability $1-\delta$, and initial MAML-RL optimality gap $\Delta_0 := J(\theta^\star) - J(\theta_0)$.

\end{thrm}
 
\noindent \textbf{Task selection bias:} We emphasize that the bias in the ergodic convergence rate comes from the task selection (steps 2-10 of Algorithm \ref{alg:metatrain}), and it can be made sufficiently small for $L = \mathcal{O}(1)$. Namely, there exists an $L$ for which that bias is minimized given arbitrarily different tasks in $\mathcal{M}$. For instance, consider the worst-case scenario where the tasks in $\mathcal{M}$ are all substantially different from each other. Then, $L = M$ guarantees that such bias is zero while recovering the convergence rate for the setting without subset selection. Hence, for the case where there are sufficiently similar tasks in $\mathcal{M}$, we let the practitioner to set $L \ll M$ and ensure that such bias remains negligible.

Let $\mathcal{S}^{\mathcal{S}}_c := \frac{L}{M}\mathcal{S}^{\mathcal{M}}_c + \mathcal{O}(Mn_s)$ and $\mathcal{S}^{\mathcal{M}}_c := \mathcal{O}(MNn_s)$ denote the total number of samples in Algorithm \ref{alg:metatrain} to find an $\epsilon$-near stationary solution, with and without task selection, respectively.

\begin{cor}(Sample complexity)\label{corollary: sample complexity local} Let the arguments of Theorem \ref{theorem:ergodic convergence} hold. Suppose the number of iterations scales as $N = \mathcal{O}(1/\epsilon)$ and the number of tasks in the task pool is sufficiently large as $M = \mathcal{O}(1/\epsilon)$. Therefore, our coreset-based task selection offers a sample complexity reduction such as $\mathcal{S}^{\mathcal{M}}_c = \textcolor{blue}{\mathcal{O}(1/\epsilon)}\mathcal{S}^{\mathcal{S}}_c$, with high probability. 
\end{cor}

\noindent \textbf{Task selection trade-off:} It is also worth highlighting the trade-off of selecting $L$ in a  heterogeneous task regime. It is evident that by setting $L$ small, when $M$ is sufficiently large, will be beneficial for reducing the number of samples when $\mathcal{O}\left(\frac{1}{M}\sum_{j \in \mathcal{M}} \min_{i \in \mathcal{S^\star}}\xi_{i,j}\right)$ is sufficiently small (i.e., order  $\epsilon$). However, when the tasks are sufficiently different, setting $L$ small may even prevent  convergence to a stationary solution. We refer the reader to \citep{mirzasoleiman2020coresets1,yang2023towards} for algorithmic alternatives that do not require a pre-specified $L$.\\

\noindent \textbf{Discussion:} Theorem \ref{theorem:ergodic convergence} and Corollary \ref{corollary: sample complexity local} summarize our main results for the MAML-RL setting. In particular, in Theorem \ref{theorem:ergodic convergence}, $\frac{1}{N}\sum_{n = 0}^{N-1}\|\nabla J(\theta_n)\|^2_2$ is controlled by two terms. The first term scales as $\mathcal{O}\left(\frac{\Delta_0}{\eta_{\text{out}}N}\right)$ and it refers to the complexity of finding a stationary solution given the initial meta-policy parameter $\theta_0$. On the other hand, as previously discussed, $\mathcal{O}\left(\frac{1}{M}\sum_{j \in \mathcal{M}} \min_{i \in \mathcal{S^\star}}\xi_{i,j}\right)$ is due to the meta-training over the weighted subset of tasks $\mathcal{S}$ instead of the entire task pool $\mathcal{M}$. 

Although \cite{mirzasoleiman2020coresets1,yang2023towards} also highlight the effect of the additive bias in the context of coresets for data-efficient deep-learning, they assume the direct computation of gradients which simplifies the setting and prevents characterization of the sample complexity and subsequently the benefit of task selection. We fill that gap for the MAML-RL problem and stress that the task selection benefit on the sample complexity reduction is not an artifact of the zeroth-order gradient estimation scheme used in this work, and it may be extended to any derivative-free approaches \citep{salimans2017evolution}.\\

\noindent \textbf{Proof idea:} The main step in the proof of Theorem \ref{theorem:ergodic convergence} is to control the gradient estimation error $\|\nabla J(\theta) - \nabla_{\mathcal{S}}J(\theta)\|$ for any $\theta \in \Theta$. To do so, we first observe that $$\|\nabla J(\theta) - \nabla_{\mathcal{S}}J(\theta)\| \leq \underbrace{\|\nabla J(\theta) - \nabla_{\mathcal{M}}J(\theta)\|}_{\text{Zeroth-order estimation error}} + \underbrace{\|\nabla J_{\mathcal{M}}(\theta) - \nabla_{\mathcal{S}}J(\theta)\|}_{\text{Task selection bias}},$$ 
where the zeroth-order estimation error can be controlled by making $n_s$ sufficiently large and $r$ sufficiently small through matrix concentration inequalities \citep{tropp2012user}. Moreover, we control the task selection bias by first using the fact that $\mathcal{F}(S) \geq (1-e^{-1})\mathcal{F}(\mathcal{S}^\star)$ \citep{nemhauser1978analysis}. Then, by also controlling the estimation error in $g_j(\theta_0)$ and $g_i(\theta_0)$ and making $C$ sufficiently small, $\|\nabla J_{\mathcal{M}}(\theta) - \nabla_{\mathcal{S}}J(\theta)\| \lesssim \epsilon + \frac{1}{M}\sum_{j \in \mathcal{M}} \min_{i \in \mathcal{S^\star}}\xi_{i,j}$, which can also be made sufficiently small by carefully tuning $L$. Subsequent proof steps follows from Assumptions \ref{assump: Lipschitz of the MAML task specific} and \ref{assump: uniform bound task specific grad}. We refer the reader to Appendix \ref{appendix: grad approx}, \ref{appendix: ergodic convergence rate} and \ref{appendix: sample complexity local} for the detailed proof. Next, we consider the MAML-LQR problem and discuss the benefit of task selection in the setting where $J_j(\theta)$ satisfies gradient dominance. 

\subsection{Linear Quadratic Regulator (LQR) Problem}

Consider the MAML-LQR problem from \cite{toso2024meta}, where the task pool $\mathcal{M}$ is composed of $M$ distinct LQR tasks $T_j = (A_j, B_j, Q_j, R_j)$ with systems matrices $A_j \in \mathbb{R}^{d_1\times d_1}$, $B_j \in \mathbb{R}^{d_1\times d_2}$, and cost matrices $Q_j \in \mathbb{S}^{d_1}_{\succeq 0}$, $R_j \in \mathbb{S}^{d_2}_{\succ 0}$, for any $j \in [M]$. In particular, each task $T_j \in \mathcal{M}$ is equipped with the objective of designing a controller $K^\star_j$ that solves
\vspace{-0.2cm}
\begin{align} \label{eq:LQR_cost}
K^\star_j =  & \argmin_{K \in \mathcal{K}_j} \left\{ \hspace{-0.1cm}J_j(K) \hspace{-0.1cm}:=\hspace{-0.1cm} \mc{E}\left[\sum_{t=0}^{\infty} x^{(j)\top}_t \left(Q_j \hspace{-0.1cm}+\hspace{-0.1cm} K^\top R_j K\right) x^{(j)}_t\right] \right\} \text{ s.t. } x^{(j)}_{t+1} \hspace{-0.1cm}=\hspace{-0.1cm} (A_j \hspace{-0.1cm}-\hspace{-0.1cm} B_jK)x^{(j)}_t, 
\end{align}
where $\mathcal{K}_j := \{K~|~ \rho(A_j - B_jK) < 1\}$ denotes the task-specific stabilizing set of controllers. 

The objective of MAML-LQR is to design $K^\star$ that stabilizes any LQR task drawn from $\mathcal{P}(\mathcal{T})$, \emph{and}, $K^\star$ should only be a few PG steps away from any unseen task-specific optimal controller. Similar to \eqref{eq:MAML_RL}, with $\theta = K$, the MAML-LQR problem is:
\vspace{-0.1cm}
\begin{align}\label{eq:MAML_LQR}
    K^\star &= \text{argmax}_{K \in \overline{\mathcal{K}}} J(K) := \mc{E}_{T_j \sim \mathcal{P}(\mathcal{T})} J_j(K - \eta_{\text{inn}} \nabla J_j(K))\notag\\
    & \text{s.t. }\quad  x^{(j)}_{t+1} = (A_j - B_jK)x^{(j)}_t, \forall T_j\sim \mathcal{P}(\mathcal{T}),
\end{align}
where $\overline{\mathcal{K}}:= \cap_{T_j\sim \mathcal{P}(\mathcal{T})} \mathcal{K}_j$ denotes the MAML-LQR stabilizing set. We note that the crucial difference between \eqref{eq:MAML_RL} and \eqref{eq:MAML_LQR} is the necessity for designing a controller $K \in \overline{\mathcal{K}}$ that \emph{stabilizes} any system $(A_j,B_j)$ drawn from the distribution of tasks $T_j \sim \mathcal{P}(\mathcal{T})$. Later, we show that Algorithm \ref{alg:metatrain} produces such stabilizing controllers, while also reducing task redundancy with task selection.

We emphasize that our goal is to understand and characterize the benefit of task selection (Algorithm \ref{alg:metatrain}), for learning stabilizing controllers that can quickly adapt to unseen tasks in the LQR setting. For this purpose, and following \cite{toso2024meta}, we next define the task specific and MAML-LQR stabilizing sub-level sets, as well as re-state the smoothness, gradient dominance and task heterogeneity properties of the LQR problem.

\vspace{-0.2cm}
\begin{definition} (Stabilizing sub-level sets) \label{def:stabilizing_set} For any task $T_j \sim \mathcal{P}(\mathcal{T})$, the task-specific sub-level set $\mathcal{G}^{\mu}_{j} \subseteq \mathcal{K}_j$ is
 $
        \mathcal{G}^{\mu}_j := \left\{K\; | \; J_j(K) - J_j(K^\star_j) \leq \mu \Delta^{(j)}_0\right\}, \text{ with } \Delta^{(j)}_0 = J^{(j)}(K_0) - J^{(j)}(K^\star_j), 
$
$\mu > 0$. In addition, the MAML-LQR stabilizing sub-level set is ${\mathcal{G}} := \cap_{j \sim \mathcal{P}(\mathcal{T})} \mathcal{G}^{\mu}_j \subseteq \overline{\mathcal{K}}$. 
\end{definition}
\vspace{-0.4cm}
\begin{assumption} (Initial stabilizing controller) \label{assumption:initial_stabilizing_K0} $K_0 \in \mathcal{G}$ \footnote{As stressed in \cite{toso2024meta},  \texttt{MAML-LQR} must be initialized from an stabilizing controller to produce finite costs and subsequently well-defined gradient estimations.}. 
\end{assumption}

\begin{assumption}\label{assump: heterogeneity} (Task heterogeneity) For any two distinct tasks $T_i,T_j \sim \mathcal{P}(\mathcal{T})$ we have that
\begin{align*}
\underset{i\neq j}{\max} \lVert  A_i -A_j\rVert \leq \epsilon_A, 
\underset{i\neq j}{\max} \lVert  B_i -B_j\rVert \leq \epsilon_B,
\underset{i\neq j}{\max} \lVert  Q_i -Q_j\rVert \leq \epsilon_Q,
\underset{i\neq j}{\max} \lVert  R_i -R_j\rVert \leq \epsilon_R,
\end{align*}
where $\epsilon_A, \epsilon_B, \epsilon_Q, \epsilon_R \geq 0$. We further denote $\epsilon_{\text{het}} = (\epsilon_A, \epsilon_B, \epsilon_Q, \epsilon_R)$.
\end{assumption}

\begin{lemma}[Lemma 4 from \cite{toso2024meta}]\label{lemma:gradient_heterogeneity} For any two distinct tasks $T_i, T_j \sim \mathcal{P}(\mathcal{T})$ and stabilizing controller $K \in \mathcal{G}$. It holds that, $\|\nabla J_i(K) - \nabla J_j(K)\| \leq {f}(\epsilon_{\text{het}})$, where ${f}(\epsilon_{\text{het}})$ denotes the gradient heterogeneity bias.
\end{lemma}

\vspace{-0.5cm}
\begin{lemma}\label{lemma:LQR} Given any task $T_j \sim \mathcal{P}(\mathcal{T})$ and stabilizing controllers $K, K^\prime \in \mathcal{G}$ such that $\|\Delta\| := \|K^{\prime} -K\|_F <\infty$. It holds that $\|\nabla J_j(K)\|_F \leq \phi$,
\begin{align*}
\left|J_j\left(K^{\prime}\right)-J_j(K)\right| \leq \beta J_j(K) \|\Delta\|_F, \quad
\left\|\nabla J_j\left(K^{\prime}\right)-\nabla J_j(K)\right\|_F \leq \psi\|\Delta\|_F,
\end{align*}
and $\|\nabla J_j(K)\|_F^2 \geq \lambda_j (J_j(K)-J_j(K^\star_j))$, where $\lambda_j > 0$ denotes the gradient dominance constant.
\end{lemma}

We remark that Lemma \ref{lemma:LQR} was initially proved in \cite{fazel2018global} and subsequently revisited in \cite{gravell2020learning, wang2023model, toso2024meta}, where the explicit expression of the problem dependent constants $\phi$, $\beta$, $\psi$ are provided.

\begin{thrm} (Gap to optimality) \label{theorem: global convergence guarantees} Suppose that Assumptions \ref{assumption: submodular function upper bound}, \ref{assumption:initial_stabilizing_K0} and \ref{assump: heterogeneity} hold. In addition, suppose that the inner and outer step-sizes are of the order $\eta_{\text{inn}} = \mathcal{O}(\epsilon/d)$ and $\eta_{\text{out}} = \mathcal{O}(1)$, and that the number of samples and smoothing radius are set according to 
\vspace{-0.2cm}
\begin{align*}
n_s \geq C_{\text{approx},2}\min(d_1,d_2)\left(\frac{\sigma^2}{\epsilon^2}+\frac{ b}{3\sqrt{\min(d_1,d_2)}\epsilon}\right) \log \left(\frac{d_1 + d_2}{\delta}\right),\quad r \leq \frac{\epsilon}{C_{\text{approx},2}\psi},
\end{align*}
for a sufficiently large universal constant $C_{\text{approx},2}$. Then, when combined with task selection the \texttt{MAML-LQR} satisfies
\vspace{-0.4cm}
\begin{align*}
J_j(K_N) - J_j(K^\star)  &\leq \left(1 -\frac{\lambda_j\eta_{\text{out}}}{4}\right)^N\Delta^{(j)}_{0} + \mathcal{O}(f^2(\epsilon_{\text{het}})),
\end{align*}
\end{thrm}
with probability $1-\delta$, for any task $T_j \sim \mathcal{P}(\mathcal{T})$ with $\Delta^{(j)}_0 = J_j(\theta^\star) - J_j(\theta_0)$.

Let $\bar{\mathcal{S}}^{\mathcal{S}}_c = \frac{L}{M}\bar{\mathcal{S}}^{\mathcal{M}}_c + \mathcal{O}(Mn_s)$ and $\bar{\mathcal{S}}^{\mathcal{M}}_c = \mathcal{O}(MNn_s)$ denote the total number of samples required in Algorithm \ref{alg:metatrain} to learn an LQR controller that is $\epsilon$-close to any task-specific optimal controller up to a heterogeneity bias, with and without task selection, respectively.
\vspace{-0.1cm}
\begin{cor} (Sample complexity) \label{corollary: sample complexity global}  Let the arguments of Theorem \ref{theorem: global convergence guarantees} hold with $L = \mathcal{O}(1)$. Suppose the number of iterations scales as $N = \mathcal{O}(\log(1/\epsilon))$ and the number of tasks in the task pool is sufficiently large as $M = \mathcal{O}(\log(1/\epsilon))$. Then, using task selection one may reduce the sample complexity to  $\bar{\mathcal{S}}^{\mathcal{M}}_c = \textcolor{blue}{\mathcal{O}(\log(1/\epsilon))}\bar{\mathcal{S}}^{\mathcal{S}}_c$, with high probability. 
\vspace{-0.2cm}
\end{cor}

The proofs are detailed in Appendix \ref{appendix: global convergence} and \ref{appendix: sample complexity global}. Note that task selection does not affect the ability of  \texttt{MAML-LQR} trained on  $\mathcal{S}$ to produce stabilizing controllers, i.e., $K_n \in \mathcal{G}$ for any iteration. We defer the stability analysis to Appendix \ref{appendix: stability analysis}. We remark that both Algorithm \ref{alg:metatrain} and \texttt{MAML-LQR} \citep[Algorithm 3]{toso2024meta} converges to a controller that is $\epsilon$-close to each task-optimal controller up to a heterogeneity bias. However, by selecting a weighted set of the most informative tasks $\mathcal{S}$, the sample complexity of learning such meta-controller is reduced by a factor of $\mathcal{O}(\log(1/\epsilon))$.

\section{Numerical Validation} \label{sec:numerics}

\begin{figure*}[ht!]
  \centering
  \begin{minipage}{.32\textwidth}
    \centering
    \includegraphics[width=1.04\textwidth]{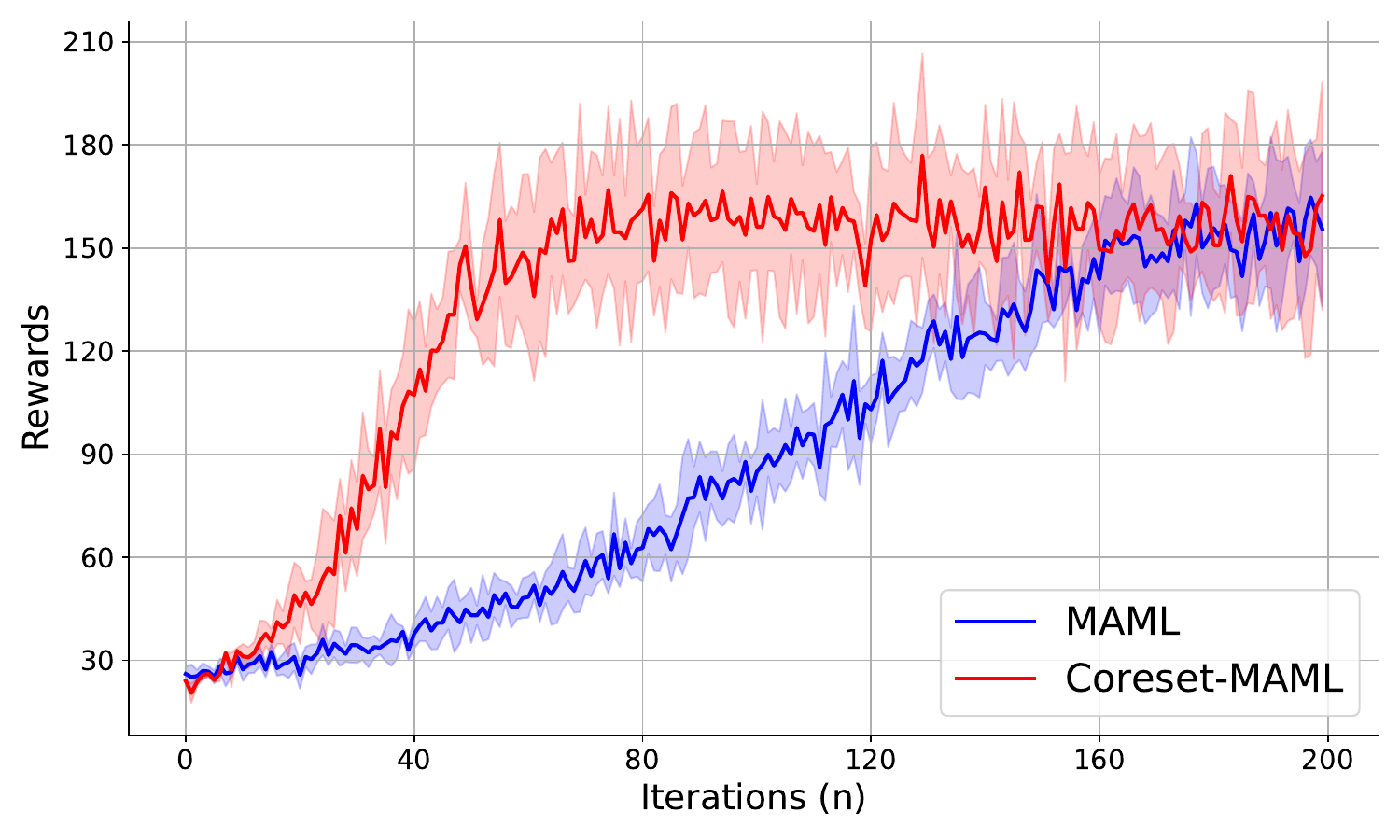}
  \end{minipage}
  \begin{minipage}{.32\textwidth}
    \centering
  \includegraphics[width=1.04\textwidth]{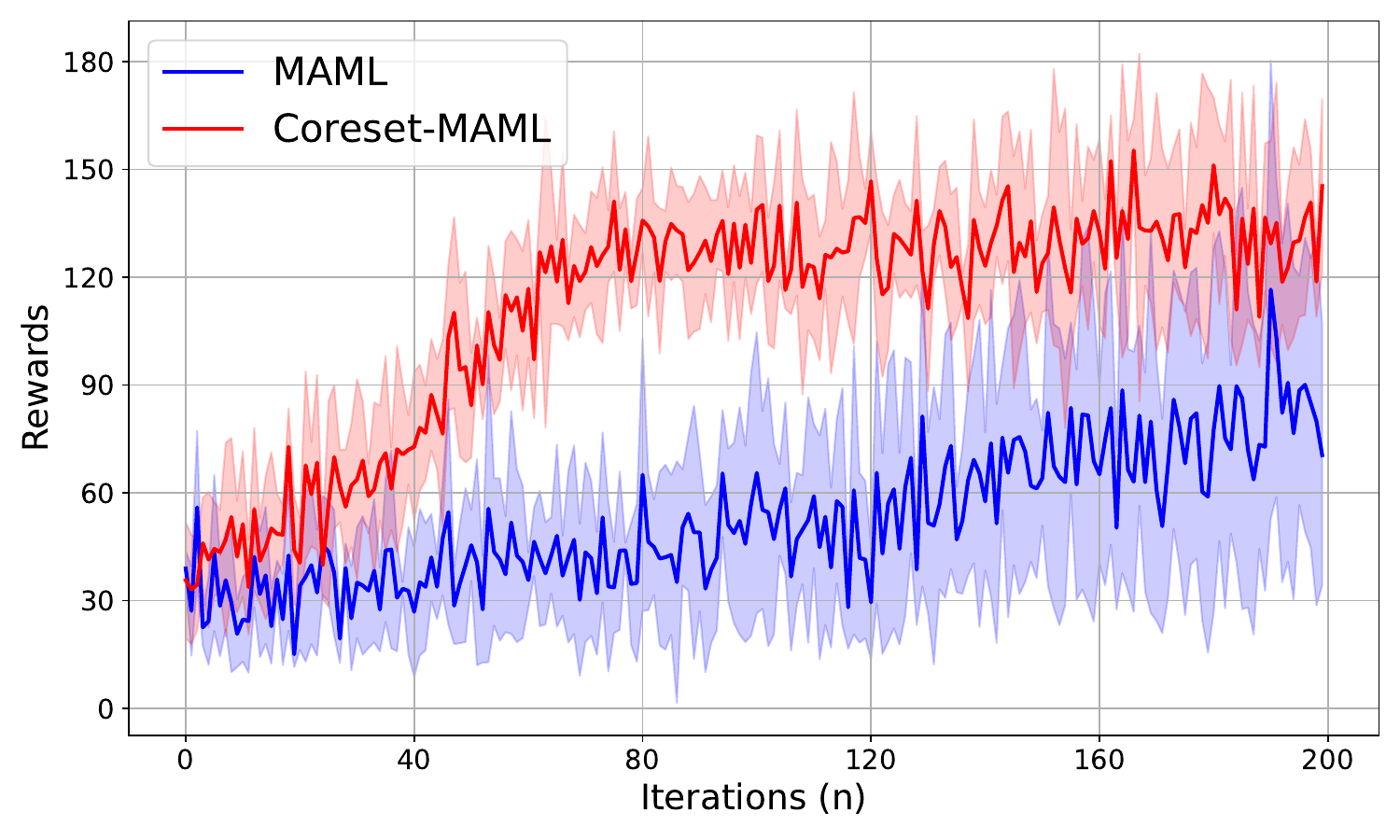}
  \end{minipage}
  \begin{minipage}{.32\textwidth}
    \centering    \includegraphics[width=1.04\textwidth]{Figures/walker.pdf}
  \end{minipage}
  \caption{Reward comparison of Algorithm \ref{alg:metatrain} and  vanilla MAML \citep{finn2017model} on Cart Pole (left), Hopper (middle), and Walker2D (right) tasks (Mujoco).}
  \label{fig:mujoco}

\end{figure*}

\noindent\textbf{Meta-Reinforcement Learning:}
We evaluate Algorithm \ref{alg:metatrain} in a deep meta-RL setting. All experiments were implemented using PyTorch, OpenAI Gym \citep{towers2024gymnasium} and Mujoco \citep{todorov2012mujoco}. In particular, we examine the cartpole, hopper, and walker2D environments \citep{towers2024gymnasium}.  These environments have physical properties of the system that vary across tasks, including the mass (cart pole, hopper, and walker2D), pole mass and length (cart pole) and friction coefficient (hopper and walker2D). The policy is parameterized by a multi-layer perceptron (MLP) architecture consisting of two hidden layers. For cart pole and hopper, we maintained an episodic task pool of 800 tasks with a selection ratio of 25\%. For each iteration, we used a batch size of 10 tasks. For walker2D, which is a more complex task with higher dimension of observation and action space, we maintained an episodic task pool of 1600 tasks with a selection ratio of 20\%. For each iteration, we used a batch size of 32 tasks. Figure \ref{fig:mujoco} shows the learning curves comparing Algorithm \ref{alg:metatrain} with vanilla MAML \cite{finn2017model}, depicting both average rewards and standard deviations across 5 runs.

Our results demonstrate that Algorithm \ref{alg:metatrain} learns approximately $4\times$ faster than the  vanilla MAML algorithm, reaching higher reward values in fewer iterations. Our approach exhibits slightly higher variance (as indicated by the shaded regions representing standard deviation) in cart pole and walker2D and it consistently outperforms the baseline in terms of learning speed across the three tasks. Notably, in hopper, the task selection method reaches a reward of 120 around iteration $60$, whereas  MAML takes approximately $200$ iterations to research a lower reward around $90$. In walker2D, the most complex task among the three tasks, our algorithm shows significant more than $4\times$ speedup of reaching higher rewards during early training stage. These empirical findings strongly support our theoretical analysis regarding sample complexity reduction (Corollary \ref{corollary: sample complexity local}), validating that careful task selection significantly enhance sample-efficiency in meta-RL. Moreover, all the experiments were running on an NVIDIA GeForce 3090 GPU with 90GB RAM. As we mentioned previously, selecting tasks via submodular maximization is efficient. In hopper, the average per-iteration running time was 3.72s for the coreset task selection algorithm and 3.35s for vanilla MAML.

\begin{figure*}[ht!]
  \centering
  \begin{minipage}{.32\textwidth}
    \centering
    \includegraphics[width=1.12\textwidth]{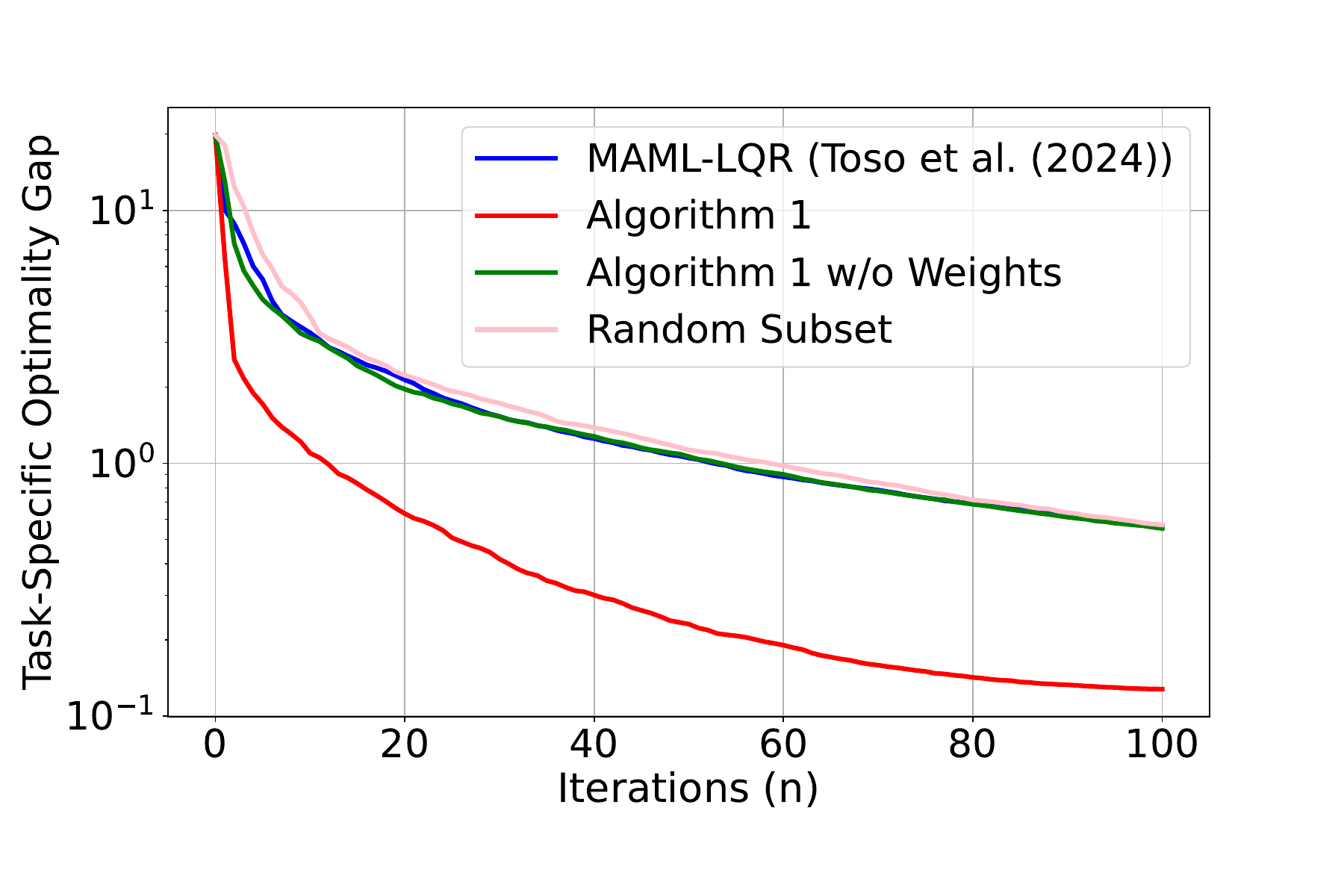}
  \end{minipage}
  \begin{minipage}{.32\textwidth}
    \centering
  \includegraphics[width=1.12\textwidth]{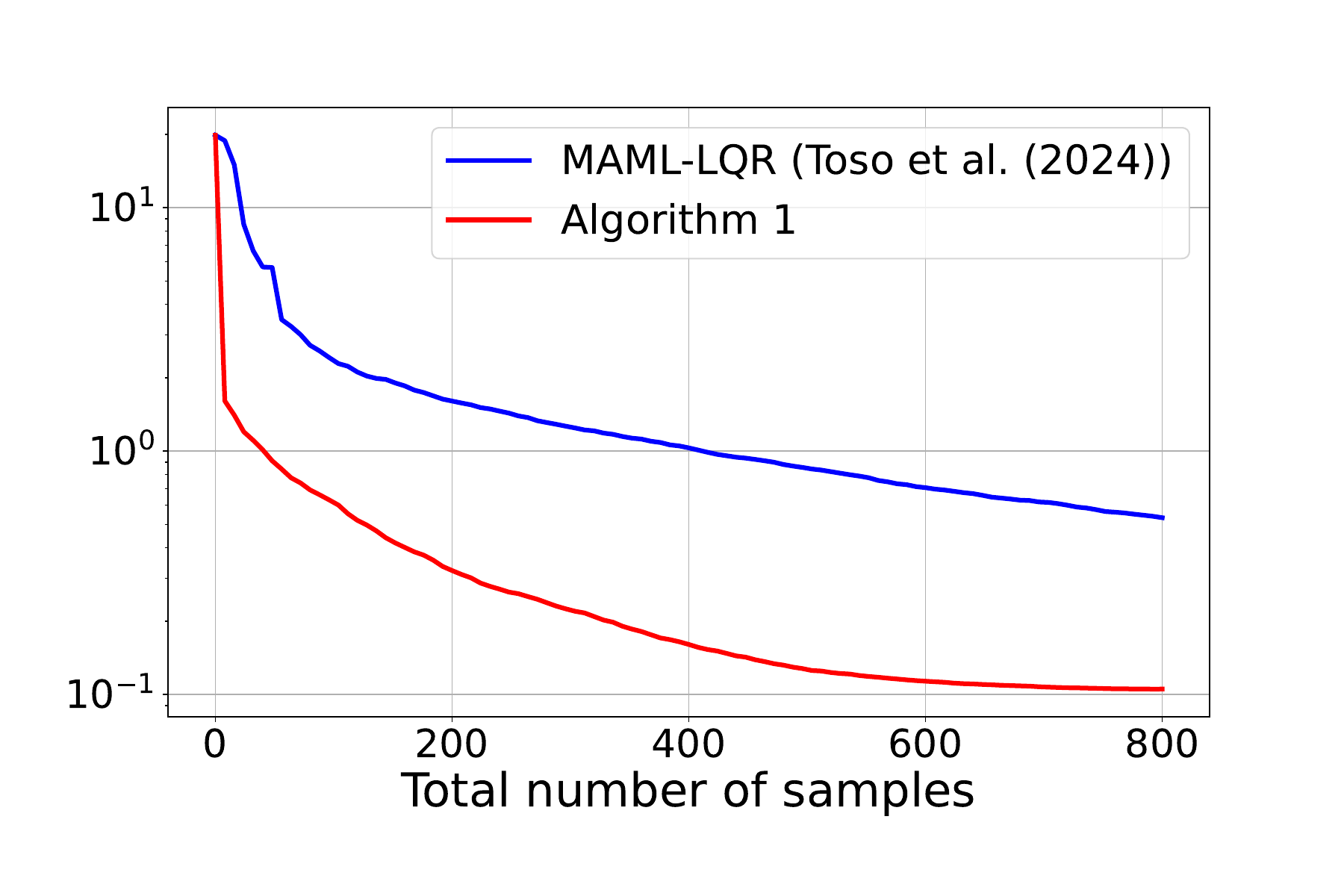}
  \end{minipage}
  \begin{minipage}{.32\textwidth}
    \centering    \includegraphics[width=1.12\textwidth]{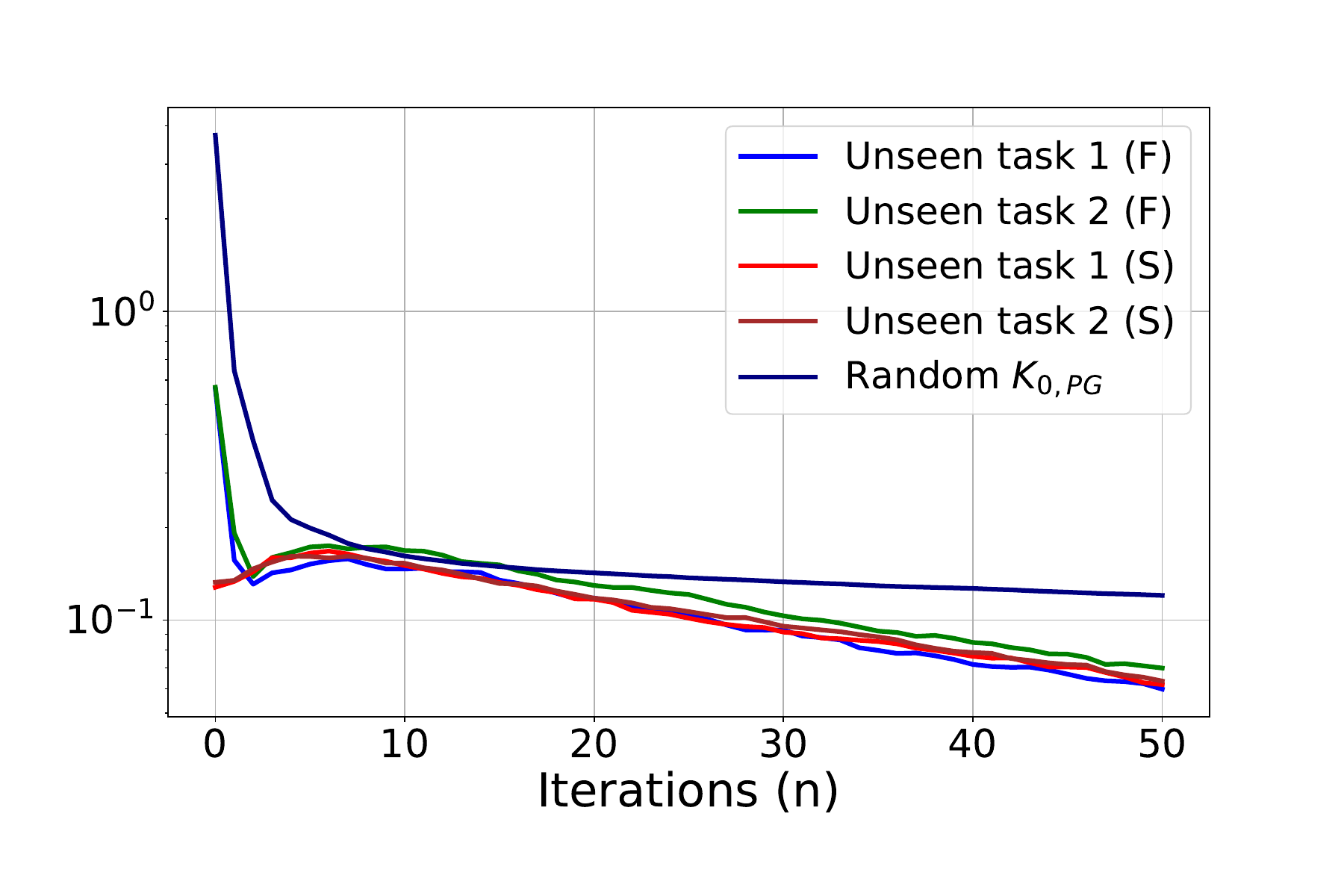}
  \end{minipage}
  \caption{Optimality gap of Algorithm \ref{alg:metatrain} in the MAML-LQR setting with respect to iterations and number of samples.}
  \label{fig:numericals}

\end{figure*}

\noindent \textbf{Linear Quadratic Regulator:} We follow the setting proposed by \citep{toso2024meta} to validate our theoretical guarantees in the MAML-LQR setting. In particular, Figure \ref{fig:numericals} (left) shows the optimality gap across iterations. We implemented the \texttt{MAML-LQR} on three scenarios: the full task pool (40 tasks), a selected subset (10 tasks), and two ablation baselines - selected subset without weight assignment and randomly selected subset. Our results demonstrate the faster convergence on the weighted selected subset, while both unweighted selected subset and random subset achieves at most the same performance as the full task pool. Moreover, Figure \ref{fig:numericals} (middle) depicts the optimality gap between the selected subset and full task pool with respect to sample size, confirming our theoretical results on sample complexity reduction for the \texttt{MAML-LQR} (Corollary \ref{corollary: sample complexity global}). Figure \ref{fig:numericals} (right) shows the task-specific optimality gap on unseen meta-testing tasks drawn from the same training task distribution $\mathcal{P(T)}$. The results show that the learned meta-controller on $\mathcal{S}$ achieves comparable generalization performance while reducing sample complexity in the meta-training, which outperforms the randomly initialized controller.

\section{Conclusions and Future Work}

We proposed a coreset-based task selection to enhance sample efficiency in meta-RL. By prioritizing the most informative and diverse tasks, Algorithm \ref{alg:metatrain} addressed the task redundancy of traditional meta-RL. We demonstrated that task selection reduces the sample complexity of finding $\epsilon$-near optimal solutions for both MAML-RL (i.e., by a factor of $\mathcal{O}(1/\epsilon)$) and MAML-LQR (i.e., proportional to $\mathcal{O}(\log(1/\epsilon))$), which are further validated through deep-RL and LQR experiments. Future work involves  adding a clustering layer, as in \cite{toso2023learning}, based on the task weights, to the meta-training pipeline to alleviate the heterogeneity bias in the MAML-LQR setting.

\section*{Acknowledgements} 

Leonardo F. Toso is funded by the Center for AI and Responsible Financial Innovation (CAIRFI) Fellowship and by the Columbia Presidential Fellowship. James Anderson is partially funded by NSF grants ECCS 2144634 and 2231350, and the Columbia Data Science Institute.

\bibliographystyle{abbrvnat}
\bibliography{bibliography}

\newpage

\section{Appendix}\label{sec:appendix}

\noindent \textbf{Roadmap:} This appendix is organized as follows: First we extend our related work section and remind the reader of the model-agnostic meta-reinforcement learning problem and the matrix Bernstein inequality from \cite{tropp2012user}, where the later is crucial for controlling the error in the MAML gradient approximation due to the zeroth-order estimation. Next, in Section \ref{appendix: grad approx}, we characterize that estimation error and prove that for a sufficiently large number of samples $n_s$, and sufficiently small smoothing radius $r$, the estimation error is composed of a sufficiently small error $\epsilon$ and an additive bias due to the task selection step in Algorithm \ref{alg:metatrain}, with high probability. Then, in Section \ref{appendix: ergodic convergence rate} we derive the ergodic convergence rate of Algorithm \ref{alg:metatrain} for the general setting of non-concave task-specific reward function $J_i(\theta)$. The sample complexity reduction benefit of task selection is then discussed in Section \ref{appendix: sample complexity local}. In Sections \ref{appendix: global convergence}, \ref{appendix: sample complexity global} and \ref{appendix: stability analysis}, we apply Algorithm \ref{alg:metatrain} to the MAML-LQR problem where $J_i(\theta)$ satisfy a gradient dominance property. %Finally in Section \ref{appendix:numerics} we provide further details on the experimental setup considered in our numerical validation section.   

\subsection{Related Work} \label{sec:related}

\begin{itemize}
\item \textbf{Meta-Reinforcement Learning:} There is a wealth of literature in meta-RL, with applications spanning robot manipulation \citep{yu2017preparing,arndt2020meta,ghadirzadeh2021bayesian}, robot locomotion \citep{song2020rapidly,yu2020learning}, build energy control \citep{luna2020information}, among others. Most relevant to our work are \cite{song2019maml,song2020rapidly} that estimate the meta-gradient through evolutionary strategy. Similarly, we consider a zeroth-order estimation of the task-specific and meta-gradients. In contrast, these works treat all the tasks equally, leading to task redundancy which we handle with a derivative-free coreset learning approach to enhance data-efficiency in meta-RL. 

\vspace{-0.2cm}
\item \textbf{Meta-Reinforcement Learning Task Selection:} Beyond the line of work on coresets for data-efficient training of machine learning models \citep{mirzasoleiman2020coresets1,killamsetty2021grad, yang2023towards, pooladzandi2022adaptive, balakrishnan2022diverse}, which use submodular optimization for subset selection, the works \citep{luna2020information,zhan2024data} are particularly relevant to this paper. In particular, \cite{zhan2024data} does not focus on RL tasks and approximates gradients using the pre-activation outputs of the last layer for classification tasks. That simplifies the problem but prevent them from deriving sample complexity guarantees. On the other hand, \cite{shin2023task} employs an information-theoretic metric to evaluate task similarities and relevance, considering a general MAML training framework rather than the policy gradient-based approach discussed here.

\vspace{-0.2cm}
\item \textbf{Model-free Learning for Control:} The linear quadratic regulator (LQR) problem has recently been taken as a fundamental baseline for establishing theoretical guarantees of policy optimization in control and reinforcement learning \citep{fazel2018global}. In particular, studies on multi-task and multi-agent learning for control \citep{zhang2023multi,wang2023model,tang2023zeroth,toso2024asynchronous,toso2024meta,lee2024nonasymptotic,lee2024regret} have derived non-asymptotic guarantees for various learning architectures within the scope of model-free LQR. Most relevant to our work are \cite{molybog2021does,musavi2023convergence,toso2024meta,aravind2024moreau,pan2024model}, which also study the meta-LQR problem and provide provable methods for learning meta-controllers that adapt quickly to unseen LQR tasks. In contrast to these works, we leverage the MAML-LQR problem as a case study to highlight the sample complexity reduction enabled by our embedded task selection approach.
\end{itemize}

\subsection{Notation and Background Results}

\noindent \textbf{Notation:}  Let $[M]$ denote the set of integers $\{1,2,\ldots,M\}$, and $\rho(\cdot)$ the spectral radius of a square matrix. Let $\|\cdot\|$ and $\|\cdot\|_F$ denote the spectral and Frobenius norm, respectively. We use $\mathcal{O}(\cdot)$ to omit constant factors in the argument. Throughout the text and when its clear from the context we use $i$ and $j$ to denote tasks $T_i$ and $T_j$. 

\vspace{0.3cm}

\noindent \textbf{Model-Agnostic Meta-Reinforcement Learning Problem:} We recall that the one-shot model-agnostic meta-reinforcement learning problem can be written as follows:

\begin{align*}
    \theta^\star := \text{argmax}_{\theta \in \Theta} J(\theta) := \mc{E}_{i \sim \mathcal{P}(\mathcal{T})} J_i(\theta + \eta_{\text{inn}} \nabla J_i(\theta)),
\end{align*}
where $J_i(\theta) := \mc{E}_{\tau \sim \mathcal{P}_i(\tau|\theta)}[\mathcal{R}_i(\tau)]$, and $\eta_{\text{inn}}$ denotes some positive step-size. 

\vspace{0.3cm}

\begin{lemma}[\cite{tropp2012user}] (Matrix Bernstein Inequality)  Let $\left\{Z_l\right\}_{l=1}^m$ be a set of $m$ independent random matrices of dimension $d = d_1 \times d_2$ with $\mathbb{E}\left[Z_l\right]=Z$, $\left\|Z_l-Z\right\| \leq b$ almost surely, and maximum variance 

$$\max \left(\left\|\mathbb{E}\left(Z_l Z_l^{\top}\right)-Z Z^{\top}\right\|,\left\|\mathbb{E}\left(Z_l^{\top} Z_l\right)-Z^{\top} Z\right\|\right) \leq \sigma^2,$$ 
and sample average $\widehat{Z}:=\frac{1}{m} \sum_{l=1}^m Z_l$. Let a small tolerance $\epsilon \geq 0$ and small probability $0 \leq \delta \leq 1$ be given. If
$$
m \geq 2\left(\frac{\sigma^2}{\epsilon^2}+\frac{ b}{3\epsilon}\right) \log \left(\frac{d_1 + d_2}{\delta}\right)
$$
$\text { then } \mathbb{P}\left[\|\widehat{Z}-Z\| \leq \epsilon\right] \geq 1-\delta \text {.}$   
\label{lemma:Bernstein}
\end{lemma}

\begin{lemma}[\textbf{Young's inequality}]  Given any two matrices $A, B \in \mathbb{R}^{d_1\times d_2}$ and a positive scalar $\nu$, it holds that
\begin{subequations}
\begin{align}\label{eq:youngs}
\|A+B\|_2^2 &\leq(1+\nu)\|A\|_2^2+\left(1+\frac{1}{\nu}\right)\|B\|_2^2 \leq (1+\nu)\|A\|_F^2+\left(1+\frac{1}{\nu}\right)\|B\|_F^2,\\
\langle A, B\rangle & \leq \frac{\nu}{2}\lVert A\rVert_2^2 +\frac{1}{2\nu}\lVert B \rVert_2^2  \leq  \frac{\nu}{2}\lVert A\rVert_F^2 +\frac{1}{2\nu}\lVert B \rVert_F^2.\label{eq:youngs_inner_product}
\end{align}
\end{subequations}
\end{lemma}

\vspace{0.3cm}

\subsection{Gradient Estimation Error} \label{appendix: grad approx}

Let $\tilde{J}(\theta):= \mc{E}_{u \sim \mathbb{S}_r}[J(\theta + u)]$ denote the Gaussian smoothing of the MAML reward function, with smoothing radius $r > 0$ and $u$ being randomly drawn from a uniform distribution over matrices of dimension $d = d_1\times d_2$ and operator norm $r$, $\mathbb{S}_r$. By using the two-point zeroth-order estimation we can define the gradient of the smoothed MAML reward function as

\begin{align*}
    \nabla \tilde{J}(\theta) := \mc{E}_{j,u}\left[\frac{d}{2r^2}\left(J_j(\theta + u + \eta_{\text{inn}} \nabla \tilde{J}_j(\theta + u)) - J_j(\theta - u + \eta_{\text{inn}}\nabla \tilde{J}_j(\theta - u))\right)u\right],
\end{align*}
where $\tilde{J}_j(\theta\pm u) = \mc{E}_{v \sim \mathbb{S}_r}[J_j(\theta \pm u + v)]$ denotes the Gaussian smoothing of the $i$-th task-specific expected reward function $J_j(\theta \pm u)$ incurred by the policy $\pi_{\theta \pm u}$. Note that the sample mean of $\nabla \tilde{J}(\theta)$ over tasks $i \sim \mathcal{P}(\mathcal{T})$ and samples $u \sim \mathbb{S}_r$ can be written as
\begin{align*}
    \widetilde{\nabla} {J}(\theta) := \frac{d}{2Mn_sr^2}\sum_{j\in \mathcal{M}} \sum_{l = 1}^{n_s}\left(J_j(\theta + u_l + \eta_{\text{inn}}\nabla \tilde{J}_j(\theta + u_l)) - J_j(\theta - u_l + \eta_{\text{inn}}\nabla \tilde{J}_j(\theta - u_l))\right)u_l.
\end{align*}

In addition, we define the two-point zeroth-order estimation of $\nabla J(\theta)$, over $\mathcal{M}$, as follows:

\begin{align*}
    \nabla_{\mathcal{M}} {J}(\theta) := \frac{d}{2Mn_sr^2}\sum_{j\in \mathcal{M}} \sum_{l = 1}^{n_s}\left(J_{j}(\theta + u_l + \eta_{\text{inn}}\widehat{\nabla} J_j(\theta)) - J_j(\theta - u_l + \eta_{\text{inn}}\widehat{\nabla} J_j(\theta))\right)u_l
\end{align*}
where 
\begin{align*}
    \widehat{\nabla} J_j(\theta) := \frac{d}{2n_sr^2}\sum_{l = 1}^{n_s}\left(J_{j}(\theta + v_l) - J_j(\theta - v_l)\right)v_l,
\end{align*}
and similarly we can define the estimation over the subset of tasks $\mathcal{S}$ as
\begin{align*}
    \nabla_{\mathcal{S}} {J}(\theta) := \frac{d}{2Mn_sr^2}\sum_{i\in \mathcal{S}} \gamma_i \sum_{l = 1}^{n_s}\left(J_{i}(\theta + u_l + \eta_{\text{inn}}\widehat{\nabla} J_i(\theta)) - J_i(\theta - u_l + \eta_{\text{inn}}\widehat{\nabla} J_i(\theta))\right)u_l,
\end{align*}
and define $g_i(\theta) = \frac{d}{2r^2n_s} \sum_{l = 1}^{n_s} \left(J_{i}(\theta + u_l + \eta_{\text{inn}}\widehat{\nabla} J_i(\theta)) - J_i(\theta - u_l + \eta_{\text{inn}}\widehat{\nabla} J_i(\theta))\right)u_l$. 

\vspace{0.3cm}

\noindent \textbf{Goal:} We aim to demonstrate that the error in approximation of the gradient over the subset of tasks $\mathcal{S}$, i.e., $\|\nabla J(\theta) - \nabla_{\mathcal{S}} J(\theta)\|$, is sufficient small if the number of samples $n_s$, the smoothing radius $r$ and the step-size $\eta_{\text{inn}}$ are set accordingly. To do so, let us first write the following 
\begin{align*}
    \|\nabla J(\theta) - \nabla_{\mathcal{S}} J(\theta)\| &= \|\nabla J(\theta) - \nabla_{\mathcal{M}} J(\theta) + \nabla_{\mathcal{M}} J(\theta) -\nabla_{\mathcal{S}}J(\theta)\|\\
    &\leq \underbrace{\|\nabla J(\theta) - \nabla_{\mathcal{M}} J(\theta)\|}_{\text{Zeroth-order gradient approximation}} +  \underbrace{\|\nabla_{\mathcal{M}} J(\theta)-\nabla_{\mathcal{S}} J(\theta)\|}_{\text{Task selection}},
\end{align*}
where we need to control the error in the gradient approximation over all the tasks in $\mathcal{M}$, and the error due to the task selection and training over the tasks in the subset $\mathcal{S}$.

\vspace{0.3cm}

\noindent $\bullet$ Zeroth-order gradient approximation $\|\nabla J(\theta) - \nabla_{\mathcal{M}} J(\theta)\|$:

\begin{align*}
    \|\nabla J(\theta) - \nabla_{\mathcal{M}} J(\theta)\| \leq \underbrace{\|\nabla J(\theta) - \nabla \tilde{J}(\theta)\|}_{\textbf{(a)}} + \underbrace{\|\nabla \tilde{J}(\theta) - \widetilde{\nabla} J(\theta)\|}_{\textbf{(b)}} +  \underbrace{\|\widetilde{\nabla} J(\theta) - \nabla_{\mathcal{M}} J(\theta)\|}_{(\textbf{c)}},
\end{align*}

\noindent $\textbf{(a)}$: 

    \begin{align*}
        \|\nabla J(\theta) - \nabla \tilde{J}(\theta)\| &= \|\mc{E}_{u \in \mathbb{S}_r}\left(\nabla J(\theta) - \nabla J(\theta + u)\right)\|\\
        &\stackrel{(i)}{\leq} \mc{E}_{j,u}\|\nabla J_j(\theta+\eta_{\text{inn}}\nabla J_j(\theta)) - \nabla J_j(\theta+u+\eta_{\text{inn}}\nabla J_j(\theta+u))\|\\
        &\stackrel{(ii)}{\leq} \mc{E}_{j,u}\psi\| \eta_{\text{inn}} \left(\nabla J_j(\theta) - \nabla J_j(\theta+u)\right) - u\|\\ 
        &\leq \psi r(1+\psi\eta_{\text{inn}}) \stackrel{(iii)}{\leq} \frac{3\psi r}{2},
    \end{align*}
where $(i)$ and $(ii)$ follows from Jensen's inequality and \eqref{eq: assump Lipschitz}, respectively. Moreover, $(iii)$ follows from selecting the step-size according to $\eta_{\text{inn}} \leq \frac{1}{2\psi}$. Therefore, by selecting the smoothing radius as $r \leq \frac{2\epsilon}{9\psi}$, we obtain $\|\nabla J(\theta) - \nabla \tilde{J}(\theta)\| \leq \frac{\epsilon}{3}$. 

\vspace{0.3cm}

 \noindent $\textbf{(b)}$: we first note that $\|\nabla \tilde{J}(\theta) - \widetilde{\nabla} J(\theta)\| = \|\mc{E}_{j,u}\widetilde{\nabla} J(\theta) - \widetilde{\nabla} J(\theta) \|$. In addition, since the tasks $T_j \in \mathcal{M}$ and samples $u \sim \mathbb{S}_r$ are drawn independently, we can use Lemma \ref{lemma:Bernstein} to control $\|\nabla \tilde{J}(\theta) - \widetilde{\nabla} J(\theta)\|$. Let us first denote $Z_l = \frac{d}{2r^2}\left(J_j(\theta + u_l + \eta_{\text{inn}}\nabla \tilde{J}_j(\theta + u_l)) - J_j(\theta - u_l + \eta_{\text{inn}}\nabla \tilde{J}_j(\theta - u_l))\right)u_l$ and $Z = \mc{E}_{j,u}[Z_l]$. 

 \begin{align*}
     \|Z_l\| &\leq  \frac{d}{2r} |J_j(\theta + u_l + \eta_{\text{inn}}\nabla \tilde{J}_j(\theta + u_l)) - J_j(\theta - u_l + \eta_{\text{inn}}\nabla \tilde{J}_j(\theta - u_l))|\\
     &\stackrel{(i)}{\leq} \frac{d\beta J_{\max}}{2r}\|2 u_l + \eta_{\text{inn}}\left(\nabla \tilde{J}_j(\theta + u_l) - \nabla \tilde{J}_j(\theta - u_l) \right)\| \\
     &\leq d\beta J_{\max} + \frac{d\beta J_{\max}\eta_{\text{inn}}}{2r}\|\mc{E}_{v \sim \mathbb{S}_r} \nabla J_j(\theta+u_l + v) - \nabla J_j(\theta - u_l + v)\|\\
     &\stackrel{(ii)}{\leq} d\beta J_{\max} + \frac{d\beta J_{\max}\eta_{\text{inn}} \psi }{2r}\|2u_j\| \leq \frac{3d\beta J_{\max}}{2},
 \end{align*}
where $(i)$ is due to \eqref{eq: assump Lipschitz} and $J_{\max} = \max_{j \in \mathcal{M}, \theta \in \Theta} J_j(\theta)$. In addition, $(ii)$ follows from the definition of the Gaussian smoothing of the task-specific reward function and from \eqref{eq: assump Lipschitz}. The last inequality is due to $\eta_{\text{inn}} \leq \frac{1}{2\psi}$. Let us now control the approximation bias $\|Z_l - Z\|$. 

\begin{align*}
    \|Z\| = \|\nabla \tilde{J}(\theta)\| \leq \frac{\epsilon}{3} + \|\nabla J(\theta)\| \leq \frac{\epsilon}{3} + \phi,
\end{align*}
which implies $\|Z_l - Z\| \leq \|Z_l\|+ \| Z\| \leq b:= \frac{3d\beta J_{\max}}{2} +   \frac{\epsilon}{3} +\phi $. We control the approximation variance $\|\mc{E}(Z_lZ_l^\top)  - ZZ^\top\|$ as follows:

\begin{align*}
    \|\mc{E}(Z_lZ_l^\top)  - ZZ^\top\| &\leq \|\mc{E}(Z_lZ_l^\top)\| + \|ZZ^\top\|\leq \max_{Z_l} (\|Z_l\|)^2 + \|Z\|^2\\
    &\leq \sigma^2 :=  \left(\frac{3d\beta J_{\max}}{2}\right)^2 + \left(\frac{\epsilon}{3} + \phi\right)^2,
\end{align*}
and by selecting $Mn_s \geq 36\left(\frac{\sigma^2}{\epsilon^2}+\frac{b}{\epsilon}\right) \log \left(\frac{d_1 + d_2}{\delta}\right)$, we have that $\|\nabla \tilde{J}(\theta) - \widetilde{\nabla} J(\theta)\|\leq \frac{\epsilon}{3}$ holds with probability $1-\delta$.

\noindent $\textbf{(c)}$:

\begin{small}
\begin{align*}
   \|\widetilde{\nabla} J(\theta) - \nabla^\prime J(\theta)\| \leq \frac{d}{2Mn_sr^2}\sum_{j\in \mathcal{M}} &\sum_{l = 1}^{n_s} \|\left(J_j(\theta + u_l + \eta_{\text{inn}}\nabla \tilde{J}_j(\theta + u_l)) - J_j(\theta + u_l + \eta_{\text{inn}}\widehat{\nabla} J_j(\theta))
   \right. \\ 
   & \left. +J_j(\theta - u_l + \eta_{\text{inn}}\widehat{\nabla} J_j(\theta)) - J_j(\theta - u_l + \eta_{\text{inn}}\nabla \tilde{J}_j(\theta - u_l))\right)u_l\|\\
   &\hspace{-2.5cm}\stackrel{(i)}{\leq} \frac{d\beta J_{\max}\eta_{\text{inn}}}{2Mn_sr}\sum_{j\in \mathcal{M}} \sum_{l = 1}^{n_s} \|\nabla \tilde{J}_j(\theta + u_l) - \widehat{\nabla} J_j(\theta)\| + \|\nabla \tilde{J}_j(\theta - u_l) - \widehat{\nabla} J_j(\theta)\|\\
   &\hspace{-2.5cm}\stackrel{(ii)}{\leq} d\beta J_{\max}\psi\eta_{\text{inn}} + \frac{d\beta J_{\max}\eta_{\text{inn}}}{Mn_sr}\sum_{j\in \mathcal{M}} \sum_{l = 1}^{n_s} \|\nabla \tilde{J}_j(\theta) - \widehat{\nabla} J_j(\theta)\|,
\end{align*}    
\end{small}
where $(i)$ is due to \eqref{eq: assump Lipschitz}, and $(ii)$ follows from adding and subtracting $\tilde{J}_j(\theta)$ and using \eqref{eq: assump Lipschitz} with the definition of Gaussian smoothing. Then, by \citep[Lemma 1]{flaxman2004online}, $\nabla \tilde{J}_j(\theta) = \mc{E}_u[ \frac{d}{r^2} J_j(\theta + u )u]$, which implies that $\nabla \tilde{J}_j(\theta) =  \mc{E} \widehat{\nabla} J_j(\theta) = \mc{E}[ \frac{d}{2r^2} \left(J_j(\theta + u ) - J_j(\theta - u )\right)u]$ due to the symmetric perturbation of the two-point zeroth-order approximation. Therefore, we proceed to control $\|\nabla \tilde{J}_j(\theta) - \widehat{\nabla} J_j(\theta)\| = \|\mc{E}_u \widehat{\nabla} J_j(\theta)  - \widehat{\nabla} J_j(\theta)\|$, by using Lemma \ref{lemma:Bernstein} as previously. In particular, we define $Z_l = \frac{d}{2r^2}\left(J_j(\theta + u_l) - J_j(\theta - u_l)\right)u_l$ and $Z = \mc{E}_{u}[Z_l]$. 

 \begin{align*}
     \|Z_k\| &\leq  \frac{d}{2r} |J_j(\theta + u_l) - J_j(\theta - u_l)| \stackrel{(i)}{\leq} \frac{d\beta J_{\max}}{2r}\|2 u_l\|
      \leq d\beta J_{\max},
 \end{align*}
where $(i)$ follows from \eqref{eq: assump Lipschitz}. In addition, we have 

\begin{align*}
\|Z\| = \|\nabla \tilde{J}_j(\theta)\| &\leq  \|\nabla \tilde{J}_j(\theta) - \nabla J_j(\theta)\| + \phi\\
&\leq \mc{E}_u\|\nabla J_j(\theta+u) - \nabla J_j(\theta)\| + \phi\\
&\leq \psi r + \phi \leq \frac{\epsilon}{3} + \phi,
\end{align*}
where the last inequality follows from $r \leq \frac{\epsilon}{3\psi}$. Then, the approximation bias and variance of the estimation in $\textbf{(c)}$ satisfy $\|Z_k - Z\| \leq b$ and $ \|\mc{E}(Z_kZ_k^\top)  - ZZ^\top\| \leq \sigma^2$, respectively. This implies that, by selecting $n_s \geq 18\left(\frac{4\sigma^2}{\epsilon^2}+\frac{b}{\epsilon}\right) \log \left(\frac{d_1 + d_2}{\delta}\right)$, $\|\nabla \tilde{J}_j(\theta) - \widehat{\nabla} J_j(\theta)\| \leq \frac{\epsilon}{6}$ holds with probability $1 - \delta$. Finally, by setting the inner step-size as $\eta_{\text{inn}} \leq \min\left\{\frac{r}{d\beta J_{\max}}, \frac{\epsilon}{6d\beta J_{\max}\psi}\right\}$, we have that $ \|\widetilde{\nabla} J(\theta) - \nabla^\prime J(\theta)\| \leq \frac{\epsilon}{3}$ holds with probability $1-\delta$. 

Therefore, by combining $\textbf{(a)}$, $\textbf{(b)}$ and $\textbf{(c)}$, and supposing that the number of samples $n_s$, the smoothing radius $r$, and step-size $\eta_{\text{inn}}$ are set as follows:
\begin{align*}
    n_s \geq 36\left(\frac{8\sigma^2}{\epsilon^2}+\frac{b}{\epsilon}\right) \log \left(\frac{d}{\delta}\right), r \leq \frac{\epsilon}{9\psi}, \text{ and } \eta_{\text{inn}} \leq \min\left\{\frac{2r}{d\beta J_{\max}}, \frac{1}{2\psi}\right\},
\end{align*}
the zeroth-order estimation error is sufficiently small, i.e., $\|\nabla J(\theta) - \nabla_{\mathcal{M}} J(\theta)\| \leq \epsilon $ holds with probability $1-2\delta$. 

It is worth noting that the number of samples, smoothing radius and inner step-size must be in the order of $n_s = \mathcal{O}(d^2/\epsilon^2)$, $r = \mathcal{O}(\epsilon)$, and $\eta_{\text{inn}} = \mathcal{O}(\epsilon/d)$, respectively, in order to ensure that the estimation error due to the zeroth-order approximation is sufficiently small, i.e, $\|\nabla J(\theta) - \nabla_{\mathcal{M}} J(\theta)\| = \mathcal{O}(\epsilon)$. We use $\mathcal{O}(\cdot)$ to omit the dependence on universal constants and only highlight the scaling the number of samples with the problem dimension $d$ and approximation error $\epsilon$.

\vspace{0.3cm}

\noindent $\bullet$ Task Selection $\|\nabla_{\mathcal{M}} J(\theta)-\nabla_{\mathcal{S}} J(\theta)\|$: To control the estimation error due to the task selection, we start by writing the following

\begin{align*}
    \|\nabla_{\mathcal{M}} J(\theta)-\nabla_{\mathcal{S}} J(\theta) \| & =  \left\| \frac{1}{M}\sum_{j \in \mathcal{M}}g_j(\theta) - \frac{1}{M}\sum_{i \in \mathcal{S}}\gamma_i g_i(\theta)\right\| = \frac{1}{M} \left\| \sum_{j \in \mathcal{M}}g_j(\theta) - \sum_{i \in \mathcal{S}}\gamma_i g_i(\theta) \right\|\\ 
    &\leq \frac{1}{M} \sum_{j \in \mathcal{M}} \min_{i \in \mathcal{S}}\left\| g_j(\theta) - g_i(\theta) \right\| \leq \frac{1}{M} \sum_{j \in \mathcal{M}} \min_{i \in \mathcal{S}} \max_{\theta \in \Theta} \left\| g_j(\theta) - g_i(\theta) \right\|.
\end{align*}

We recall that the greedy subset selection (i.e., steps 2-10 in Algorithm \ref{alg:metatrain}) returns a subset $\mathcal{S}$  that is a suboptimal solution of the following submodular maximization
\begin{align*}
\mathcal{S}^{\star} =\underset{\mathcal{S} \subseteq \mathcal{M}}{\text{argmax} } \;\ \mathcal{F}(\mathcal{S}):=C - \sum_{j \in \mathcal{M}} \min _{i \in \mathcal{S}}\left\| {g}_j(\theta_0) -  g_i(\theta_0)\right\|, \text { subject to }|\mathcal{S}| \leq L,
\end{align*}
for any $\theta_0 \in \Theta$. In particular, by \citep[Section 4]{nemhauser1978analysis}, we know that the value of the greedy optimization is close to the optimal as $\mathcal{F}(\mathcal{S}) \geq (1 - e^{-1})\mathcal{F}(\mathcal{S^\star})$. This implies that
\begin{align*}
    \sum_{j \in \mathcal{M}} \min _{i \in \mathcal{S}}\left\| {g}_j(\theta_0) -  g_i(\theta_0)\right\| \leq Ce^{-1} + (1 - e^{-1})\sum_{j \in \mathcal{M}} \min _{i \in \mathcal{S^\star}}\left\| {g}_j(\theta_0) -  g_i(\theta_0)\right\|, 
\end{align*}
and by taking the maximum of both sides with respect to $\theta \in \Theta$, we obtain
\begin{align*}
    \frac{1}{M}\sum_{j \in \mathcal{M}} \min _{i \in \mathcal{S}} \max_{\theta \in \Theta} \left\| {g}_j(\theta) -  g_i(\theta)\right\| \leq \frac{Ce^{-1}}{M} + \frac{(1 - e^{-1})}{M}\sum_{j \in \mathcal{M}} \min _{i \in \mathcal{S^\star}} \max_{\theta \in \Theta}\left\| {g}_j(\theta) -  g_i(\theta)\right\|,
\end{align*}
where we control $\left\| {g}_j(\theta) -  g_i(\theta)\right\|$ as follows

\begin{align*}
   \left\| {g}_j(\theta) -  g_i(\theta)\right\| &\leq   \left\| {g}_j(\theta) -  \nabla J_j(\theta +\eta_{\text{inn}}\nabla J_j(\theta))\right\| + \left\| \nabla J_j(\theta +\eta_{\text{inn}}\nabla J_j(\theta)) - \nabla J_j(\theta)\right\|\\
   &+\left\| {g}_i(\theta) -  \nabla J_i(\theta +\eta_{\text{inn}}\nabla J_i(\theta))\right\| + \left\| \nabla J_i(\theta +\eta_{\text{inn}}\nabla J_i(\theta)) - \nabla J_i(\theta)\right\|\\
   &+ \|\nabla J_i(\theta) - \nabla J_j(\theta)\| \\
   &\stackrel{(i)}{\leq}   \left\| {g}_j(\theta) -  \nabla J_j(\theta +\eta_{\text{inn}}\nabla J_j(\theta))\right\| +\left\| {g}_i(\theta) -  \nabla J_i(\theta +\eta_{\text{inn}}\nabla J_i(\theta))\right\|\\
   &+ 2\eta_{\text{inn}}\psi \phi + \|\nabla J_i(\theta) - \nabla J_j(\theta)\|\\
   &\stackrel{(ii)}{\leq} \left\| {g}_j(\theta) -  \nabla J_j(\theta +\eta_{\text{inn}}\widehat{\nabla} J_j(\theta))\right\| +\left\| {g}_i(\theta) -  \nabla J_i(\theta +\eta_{\text{inn}}\widehat{\nabla} J_i(\theta))\right\|\\
   &+ \eta_{\text{inn}}\psi \|\widehat{\nabla} J_j(\theta) - \nabla J_j(\theta)\| + \eta_{\text{inn}}\psi \|\widehat{\nabla} J_i(\theta) - \nabla J_i(\theta)\|\\
   &+  2\eta_{\text{inn}}\psi \phi + \|\nabla J_i(\theta) - \nabla J_j(\theta)\|,
\end{align*}
where $(i)$ and $(ii)$ are due to \eqref{assump: Lipschitz of the MAML task specific}. Therefore, we note that for either inner and outer zeroth-order gradient approximations, i.e., $\left\| {g}_j(\theta) -  \nabla J_j(\theta +\eta_{\text{inn}}\widehat{\nabla} J_j(\theta))\right\|$ and $\|\widehat{\nabla} J_j(\theta) - \nabla J_j(\theta)\|$, respectively, the number of samples and smoothing radius can be set according to 
\begin{align*}
n_s \geq C_{\text{approx},1}\left(\frac{\sigma^2}{\epsilon^2}+\frac{ b}{3\epsilon}\right) \log \left(\frac{d_1 + d_2}{\delta}\right), r \leq \frac{\epsilon}{C_{\text{approx},1}\psi}
\end{align*}
to guarantee a small estimation error of $\epsilon$, with probability $1-\delta$, where $C_{\text{approx},1}$ is some positive universal constant, and $\sigma^2 =  \left(d\beta J_{\max}\right)^2 + \left(\epsilon + \phi\right)^2$, $b = d\beta J_{\max} +   \epsilon +\phi$. Therefore, following Definition \ref{definition: maximum normed difference}, i.e., $\xi_{i,j} := \max_{\theta \in \Theta}\|\nabla J_i(\theta) - \nabla J_j(\theta)\|$, the estimation error due to the task selection is bounded as follows:
\begin{align*}
    \|\nabla_{\mathcal{M}} J(\theta)-\nabla_{\mathcal{M}} J(\theta)\| \leq \epsilon + \frac{1}{M}\sum_{j \in \mathcal{M}} \min_{i \in \mathcal{S^\star}}\xi_{i,j},
\end{align*}
which holds with high probability $1 - \delta$ for $n_s$ and $r$ as above, and $C \leq \frac{Me\epsilon }{C_{\text{approx},1}}$. Finally, the gradient estimation error $\|\nabla J(\theta) - \nabla_{\mathcal{S}} J(\theta)\|$ in Algorithm \ref{alg:metatrain}, is controlled by a sufficiently small error that comes from the zeroth-order estimation and an additive bias due to the task selection. That is,
\begin{align*}
    \|\nabla J(\theta) - \nabla_{\mathcal{S}} J(\theta)\| \leq \epsilon + \frac{1}{M}\sum_{j \in \mathcal{M}} \min_{i \in \mathcal{S^\star}}\xi_{i,j},
\end{align*}
which holds with high probability for $n_s = \mathcal{O}(d^2/\epsilon^2)$, $r = \mathcal{O}(\epsilon)$, $\eta_{\text{inn}} = \mathcal{O}(\epsilon/d)$, and $C = \mathcal{O}(\epsilon)$. We emphasize that setting $C$ sufficiently small, i.e., $C = \mathcal{O}(\epsilon)$ is standard in the literature of coresets for data-efficient machine-learning \citep{mirzasoleiman2020coresets1, pooladzandi2022adaptive, yang2023towards} and it guarantees that the gradient estimation error due to the subset selection is sufficiently small. 

\begin{remark}(Expected Reward vs Empirical Reward) It is worth noting that our previous derivations assume that we have access to an oracle that provides the task-specific expected reward $J_j(\theta)$ incurred by any policy $\pi_{\theta}$, with $\theta \in \Theta$. However, in practice, we often do not have access to the true distribution of trajectories conditioned on the policy, i.e., $\tau \sim \mathcal{P}_i(\tau|\theta)$, that is needed to compute $J_j(\theta)$. Therefore, one may approximate the expected reward $J_j(\theta)$ with an empirical reward $\hat{J}_j(\theta) := \frac{1}{n_\tau}\sum_{l=1}^{n_\tau}[\mathcal{R}_j(\tau_l)]$, where $\{\tau_l\}^{n_\tau}_{l=1}$ are the trajectories obtained by playing with $\pi_{\theta}$, $n_\tau$ times. Note that, we can control the error between $J_j(\theta)$ and $\hat{J}_j(\theta)$ with $n_\tau$. Then, since that error should enters the analysis of Algorithm \ref{alg:metatrain} for either with or without task selection settings, we may assume the access to $J_j(\theta)$, for simplicity, but we stress that our results can be readily extended to the practical setting of empirical rewards by controlling such error with a sufficiently large $n_\tau$. 
\end{remark}

\subsection{Proof of Theorem \ref{theorem:ergodic convergence} (Ergodic Convergence Rate)} \label{appendix: ergodic convergence rate}

Recall that the meta-policy parameter is updated as follows:

\begin{align*}
    \theta_{n+1} = \theta_{n} +  \eta_{\text{out}}\nabla_{\mathcal{S}}J(\theta_n) =\theta_{n} +  \frac{\eta_{\text{out}}}{M}\sum_{i \in \mathcal{S}} \gamma_i g_i(\theta_n).  
\end{align*}

In addition, by using the definition of the meta-gradient and the gradient Lipschitz assumption \eqref{eq: assump Lipschitz}, we have that

\begin{align*}
   \|\nabla J(\theta_1) - \nabla J(\theta_2)\| &= \|\mc{E}_{i}\nabla J_i(\theta_1 + \eta_{\text{inn}}\nabla J_i(\theta_1)) - \mc{E}_{i}\nabla J_i(\theta_2 + \eta_{\text{inn}}\nabla J_i(\theta_2))\|\\
   &\stackrel{(i)}{\leq} \mc{E}_{i}\psi \|(\theta_1 - \theta_2) + \eta_{\text{inn}}\left(\nabla J_i(\theta_1) - \nabla J_i(\theta_2)\right)\|\\
   &\leq \psi \|\theta_1 - \theta_2\| + \eta_{\text{inn}}\|\nabla J_i(\theta_1) - \nabla J_i(\theta_2)\|\\
   &\stackrel{(ii)}{\leq} \psi(1+\eta_{\text{inn}}\psi) \|\theta_1 - \theta_2\|\\
   &\stackrel{(iii)}{\leq} \frac{3\psi}{2}\|\theta_1 - \theta_2\| = \overline{\psi}\|\theta_1 - \theta_2\|,
\end{align*}
for any $\theta_1, \theta_2 \in \Theta$. Here, $(i)$ and $(ii)$ follows from \eqref{assump: uniform bound task specific grad}, and $(iii)$ is due to $\eta_{\text{inn}} \leq \frac{1}{4\psi}$. Therefore, the MAML reward function incurred by policy $\pi_{\theta}$ is $\overline{\psi}$-smooth and satisfy
\begin{align*}
    J(\theta_{n}) - J(\theta_{n+1}) &\leq \langle  \nabla J(\theta_n), \theta_{n} - \theta_{n+1}\rangle + \frac{\overline{\psi}}{2}\|\theta_{n+1} -\theta_n\|^2\\
    &= \langle  \nabla J(\theta_n),  -\eta_{\text{out}}\nabla_{\mathcal{S}}J(\theta_n) + \eta_{\text{out}}{\nabla}J(\theta_n) - \eta_{\text{out}}{\nabla}J(\theta_n)  \rangle + \frac{\eta_{\text{out}}^2\overline{\psi}}{2}\|\nabla_{\mathcal{S}}J(\theta_n)\|^2\\
    & =  - \eta_{\text{out}}\|\nabla J(\theta_n)\|^2 +  \eta_{\text{out}}\langle  \nabla J(\theta_n),  \nabla J(\theta_n) -\nabla_{\mathcal{S}}J(\theta_n) \rangle + \frac{\eta_{\text{out}}^2\overline{\psi}}{2}\|\nabla_{\mathcal{S}}J(\theta_n)\|^2\\
    &\stackrel{(i)}{\leq} - \frac{\eta_{\text{out}}}{2}\|\nabla J(\theta_n)\|^2 +  \frac{\eta_{\text{out}}}{2}\|\nabla J(\theta_n) -\nabla_{\mathcal{S}}J(\theta_n)\|^2 + \frac{\eta_{\text{out}}^2\overline{\psi}}{2}\|\nabla_{\mathcal{S}}J(\theta_n)\|^2\\
    &\stackrel{(ii)}{\leq}  - \frac{\eta_{\text{out}}}{2}\|\nabla J(\theta_n)\|^2 +  \frac{\eta_{\text{out}}}{2}\|\nabla J(\theta_n) -\nabla_{\mathcal{S}}J(\theta_n)\|^2\\
    &+ \eta_{\text{out}}^2\overline{\psi}\left(\|\nabla J(\theta_n)\|^2 + \|\nabla J(\theta_n) - \nabla_{\mathcal{S}}J(\theta_n)\|^2\right),
\end{align*}
where $(i)$ and $(ii)$ follows from Young's inequalities \eqref{eq:youngs_inner_product} and \eqref{eq:youngs}, respectively. Then, by setting $\eta_{\text{out}} \leq \frac{1}{2\overline{\psi}}$ and re-arranging the terms, we can write

\begin{align*}
    \frac{\eta_{\text{out}}}{4}\|\nabla J(\theta_n)\|^2_2 &\leq J(\theta_{n+1}) - J(\theta_{n}) + \frac{3\eta_{\text{out}}}{4}\|\nabla J(\theta_n) - \nabla_{\mathcal{S}} J(\theta_n)\|^2,
\end{align*}
which can be unrolled over the iterations $n = \{0,1,\ldots, N-1\}$ to obtain 
\begin{align*}
    \frac{1}{N}\sum_{n = 0}^{N-1}\|\nabla J(\theta_n)\|^2_2 &\leq   \frac{4\Delta_0}{\eta_{\text{out}}N}  +  6\left(\frac{1}{M}\sum_{j \in \mathcal{M}} \min_{i \in \mathcal{S^\star}}\xi_{i,j}\right)^2,
\end{align*}
where we also use the fact that $J(\theta^\star) \geq J(\theta_N)$ above. In addition, we disregard $\mathcal{O}(\epsilon^2)$ since it is negligible for small $\epsilon$. We also note that the first term denotes the local algorithm's complexity to find an stationary solution given an initial optimality gap $\Delta_0 = J(\theta^\star) - J(\theta_{0})$, and the second term scales with the gradient approximation error due to the zeroth-order estimation and task selection. Therefore, by setting the number of iterations as $N = \mathcal{O}(1/\epsilon)$, Algorithm \ref{alg:metatrain} satisfy
\begin{align}\label{eq:stationary_solution}
    \frac{1}{N}\sum_{n = 0}^{N-1}\|\nabla J(\theta_n)\|^2_2 &\leq  \mathcal{O}\left(\epsilon +\left(\frac{1}{M}\sum_{j \in \mathcal{M}} \min_{i \in \mathcal{S^\star}}\xi_{i,j}\right)^2\right).
\end{align}

\subsection{Proof of Corollary \ref{corollary: sample complexity local} (Sample Complexity)}  \label{appendix: sample complexity local}

We let $\mathcal{S}^{\mathcal{S}}_c$ and $\mathcal{S}^{\mathcal{M}}_c$ denote the total number of samples required in Algorithm \ref{alg:metatrain} to find an $\epsilon$-near stationary solution, with and without task selection, respectively. In particular, $\mathcal{S}^{\mathcal{S}}_c = \mathcal{O}\left(LNn_s\right) + \mathcal{O}\left(Mn_s\right)$ and $\mathcal{S}^{\mathcal{M}}_c = \mathcal{O}\left(MNn_s\right)$. Note that in order to guarantee \eqref{eq:stationary_solution} with high probability, we need $n_s = \mathcal{O}(d^2/\epsilon^2)$ samples in the zeroth-order gradient estimation. Therefore, we have that $\mathcal{S}^{\mathcal{M}}_c = \mathcal{O}\left(MNn_s\right) = \mathcal{O}(d^2/\epsilon^{3})$, whereas $\mathcal{S}^{\mathcal{S}}_c = \frac{L}{M} \mathcal{O}\left(d^2/\epsilon^3\right) + \mathcal{O}\left(d^2/\epsilon^2\right).$
Then, for a sufficiently large amount of tasks in the task pool $\mathcal{M}$ (e.g., scaling as $M = \mathcal{O}(1/\epsilon)$) and a sufficiently small number of informative tasks in the subset $\mathcal{S}$, i.e., $L = \mathcal{O}(1)$, the task selection benefits from a sample complexity reduction by a factor of $\textcolor{blue}{\mathcal{O}(1/\epsilon)}$ when compared to the setting without task selection.

\subsection{Proof of Theorem \ref{theorem: global convergence guarantees} (Optimality Gap)} \label{appendix: global convergence} 

We first note Lemma \ref{lemma:LQR} is provided in the Frobenius norm. Therefore, we use the fact that for any matrix $A \in \mathbb{R}^{d_1\times d_2}$, $\|A\|_F \leq \sqrt{\min (d_1,d_2)}\|A\|$, to adapt the gradient approximation error $\|\nabla J(\theta) - \nabla_{\mathcal{S}}J(\theta)\|$ as discussed previously in the RL setting, for the MAML-LQR where $\theta = K$ and $u \sim \mathbb{S}_r$ has Frobenius norm $r$ in \texttt{ZO2P}($\cdot$). Moreover, we recall that the meta-controller is updated as follows:

\begin{align*}
    K_{n+1} = K_{n} -  \eta_{\text{out}}\nabla_{\mathcal{S}}J(K_n) =\theta_{n} -  \frac{\eta_{\text{out}}}{M}\sum_{i \in \mathcal{S}} \gamma_i g_i(K_n),  
\end{align*}
where $g_i(K_n) = \frac{d}{2r^2n_s} \sum_{l = 1}^{n_s} \left(J_{i}(K_n + u_l - \eta_{\text{inn}}\widehat{\nabla} J_i(K_n)) - J_i(K_n - u_l - \eta_{\text{inn}}\widehat{\nabla} J_i(K_n))\right)u_l$ for any task $T_i \sim \mathcal{P}(\mathcal{T})$. Then, by the gradient smoothness property in Lemma \ref{lemma:LQR}, we can write

\begin{align*}
    J_j(K_{n+1}) &- J_j(K_n) \leq \langle \nabla J_j(K_n), K_{n+1} - K_{n} \rangle + \frac{\psi}{2}\|K_{n+1} - K_n\|^2_F \\
    &=\langle \nabla J_j(K_n), -\eta_{\text{out}} \nabla_{\mathcal{S}}J(K_n) - \eta_{\text{out}} \nabla J_j(K_n) + \eta_{\text{out}} \nabla J_j(K_n) \rangle + \frac{\psi\eta_{\text{out}}^2}{2}\|\nabla_{\mathcal{S}}J(K_n)\|^2_F\\
    & \leq -\frac{\eta_{\text{out}}}{2}\|\nabla J_j(K_n)\|^2_F + \frac{\eta}{2}\|\nabla_{\mathcal{S}}J(K_n) - \nabla J_j(K_n)\|^2_F + \frac{\psi \eta_{\text{out}}^2}{2}\|\nabla_{\mathcal{S}}J(K_n)\|^2_F\\
    &\leq -\frac{\eta_{\text{out}}}{4}\|{\nabla} J_j(K_n)\|^2_F + \frac{3\eta_{\text{out}}}{4}\|\nabla_{\mathcal{S}}J(K_n) - \nabla J_j(K_n)\|^2_F,
\end{align*}
where the last two inequalities are due to Young's inequality \eqref{eq:youngs_inner_product} and \eqref{eq:youngs}, and $\eta_{\text{out}} \leq \frac{1}{4\psi}$. Let us now proceed to control the error in the meta-gradient approximation over $\mathcal{S}$, i.e.,  $\nabla_{\mathcal{S}}J(K_n)$, with respect to the task-specific gradient $\nabla J_j(K_n)$.

\begin{align*}
    \|\nabla_{\mathcal{S}}J(K_n) - \nabla J_j(K_n)\|_F &\leq \underbrace{\left\|\nabla_{\mathcal{S}}J(K_n) - \nabla J(K_n)\right\|_F}_{\text{Gradient approximation}} + \underbrace{\left\|\nabla J(K_n) - \nabla J_j(K_n)\right\|_F}_{\text{Task heterogeneity}},
\end{align*}

$\bullet$ \textbf{Task heterogeneity:} 

\begin{align*}
    &\left\|\nabla J(K_n) - \nabla J_j(K_n)\right\|_F = \left\|\mc{E}_i\nabla J_i(K_n - \eta_{\text{inn}}\nabla J_i(K_n)) - \nabla J_j(K_n)\right\|_F \\
    &\leq \mc{E}_i\left\|\nabla J_i(K_n - \eta_{\text{inn}}\nabla J_i(K_n)) - \nabla J_i(K_n)  + \nabla J_i(K_n) - \nabla J_j(K_n)\right\|_F\\
    &\leq \mc{E}_i\left\|\nabla J_i(K_n - \eta_{\text{inn}}\nabla J_i(K_n)) - \nabla J_i(K_n)\|_F  + \mc{E}_i\| \nabla J_i(K_n) - \nabla J_j(K_n)\right\|_F\\
    &\stackrel{(i)}{\leq}  \eta_{\text{inn}}\psi\mc{E}_i\|\nabla J_i(K_n)\|  + \mc{E}_i\| \nabla J_i(K_n) - \nabla J_j(K_n)\|_F\\
    &\stackrel{(ii)}{\leq} \eta_{\text{inn}}\psi\phi  + f(\epsilon_{\text{het}}),
\end{align*}
where $(i)$ is due to Lemma \ref{lemma:LQR} and $(ii)$ follows from Lemmas \ref{lemma:gradient_heterogeneity} and \ref{lemma:LQR}. 

\vspace{0.3cm}
$\bullet$ \textbf{Gradient approximation:} 

\begin{align*}
    \left\|\nabla_{\mathcal{S}}J(K_n) - \nabla J(K_n)\right\|_F  &=  \left\|\nabla_{\mathcal{S}}J(K_n) - \nabla_{\mathcal{M}} J (K_n) + \nabla_{\mathcal{M}} J (K_n) -  \nabla J(K_n)\right\|_F\\
    &\leq \underbrace{\left\|\nabla_{\mathcal{S}}J(K_n) - \nabla_{\mathcal{M}} J (K_n)\right\|_F}_{\text{Task selection}} + \underbrace{\|\nabla_{\mathcal{M}} J (K_n) -  \nabla J(K_n)\|_F}_{\text{Zeroth-order approximation}}.
\end{align*}

Following the previous analysis for the gradient estimation error of Algorithm \ref{alg:metatrain} in Section \ref{appendix: grad approx}, we know that $\left\|\nabla_{\mathcal{S}}J(K_n) - \nabla J(K_n)\right\|_F$ satisfy
\begin{align*}
    \left\|\nabla_{\mathcal{S}}J(K_n) - \nabla J(K_n)\right\|_F \leq   \epsilon + \frac{1}{M}\sum_{j \in \mathcal{M}} \min_{i \in \mathcal{S^\star}}\xi_{i,j},
\end{align*}
with probability $1 - \delta$, if the number of samples and smoothing radius are selected according to
\begin{align}\label{eq: ns and r}
    n_s \geq C_{\text{approx},2}\min(d_1,d_2)\left(\frac{\sigma^2}{\epsilon^2}+\frac{ b}{3\sqrt{\min(d_1,d_2)}\epsilon}\right) \log \left(\frac{d_1 + d_2}{\delta}\right), r \leq \frac{\epsilon}{C_{\text{approx},2}\psi},
\end{align}
with $\eta_{\text{inn}} = \mathcal{O}(\epsilon)$ and $C = \mathcal{O}(\epsilon)$. We also note that $\min(d_1,d_2)$ and $\sqrt{\min(d_1,d_2)}$ come from the Frobenius norm in Lemma \ref{lemma:LQR}. Then, we can write 
\begin{align*}
    \hspace{-0.1cm}J_j(K_{n+1}) &- J_j(K_n)  \hspace{-0.1cm}\leq\hspace{-0.1cm} -\frac{\eta_{\text{out}}}{4}\|{\nabla} J_j(K_n)\|^2_F + 3\eta_{\text{out}}\left(\eta_{\text{inn}}^2\psi^2\phi^2 \hspace{-0.1cm} + \hspace{-0.1cm} f^2(\epsilon_{\text{het}}) \hspace{-0.1cm} + \hspace{-0.1cm}\epsilon^2 \hspace{-0.1cm}+\hspace{-0.1cm} \left(\frac{1}{M}\sum_{j \in \mathcal{M}} \min_{i \in \mathcal{S^\star}}\xi_{i,j}\right)^2\right)\\
    &\stackrel{(i)}{\leq} -\frac{\lambda_j\eta_{\text{out}}}{4}\left(J_j(K_n) - J_j(K^\star_j)\right) + 3\eta_{\text{out}}\left(f^2(\epsilon_{\text{het}}) + \left(\frac{1}{M}\sum_{j \in \mathcal{M}} \min_{i \in \mathcal{S^\star}}\xi_{i,j}\right)^2 \right)\\
    &\stackrel{(ii)}{\leq} \frac{\lambda_j\eta_{\text{out}}}{4}\left(J_j(K_n) - J_j(K^\star_j)\right) + 6\eta_{\text{out}}f^2(\epsilon_{\text{het}}),
\end{align*}
where $(i)$ follows from the gradient dominance property in Lemma \ref{lemma:LQR}, $\eta_{\text{inn}} = \mathcal{O}(\epsilon)$ and disregarding $\mathcal{O}(\epsilon^2)$ since it is negligible for small $\epsilon$. $(ii)$ is due to the fact that $\frac{1}{M} \sum_{j \in \mathcal{M}}\min_{i \in \mathcal{S}^\star}\xi_{i,j} \leq f(\epsilon_{\text{het}})$.  Then, we can add and subtract $J_j(K^\star_j)$ on the LHS to obtain 
\begin{align*}
    \Delta^{(j)}_{n+1}  \leq \left(1 -\frac{\lambda_j\eta_{\text{out}}}{4}\right)\Delta^{(j)}_{n} +6\eta_{\text{out}}f^2(\epsilon_{\text{het}}),
\end{align*}
where $\Delta^{(j)}_{n} = J_j(K_n) - J_j(K^\star_j)$. Therefore, by unrolling the above expression over the iterations, $n = \{0,1,\ldots, N-1\}$, we obtain 
\begin{align}\label{eq:global convergence analysis}
    \Delta^{(j)}_{N}  &\leq \left(1 -\frac{\lambda_j\eta_{\text{out}}}{4}\right)^N\Delta^{(j)}_{0} + \frac{24}{\lambda_j}f^2(\epsilon_{\text{het}})\leq \epsilon + \frac{24}{\lambda_j}f^2(\epsilon_{\text{het}}),
\end{align}
where the last inequality is due to $N \geq \frac{4}{\eta_{\text{out}} \lambda_j}\log\left(\frac{\Delta^{(j)}_0}{\epsilon}\right)$. Then, we can conclude that Algorithm \ref{alg:metatrain}, for the MAML-LQR problem, learns a meta-controller $K_N$ that is $\epsilon$-close to any task-specific optimal controller (i.e.,  $ T_j \sim \mathcal{P}(\mathcal{T})$) up to a task heterogeneity bias that scales as $\mathcal{O}\left(f^2(\epsilon_{\text{het}})\right)$.

%\begin{remark} It is worth noting that if we set $L = M$ in Algorithm \ref{alg:metatrain}, we recover the results for MAML-LQR without task selection \citep[Theorem 4]{toso2024meta}. On the other hand, if $\frac{L}{M}$ is set accordingly to guarantee $\frac{L}{M} \geq 1 - \epsilon^\prime$, then Algorithm \ref{alg:metatrain} also benefits from a bias reduction in the order of $\mathcal{O}\left(\frac{L^2}{M^2}\right) = \mathcal{O}((1-\epsilon^\prime)^2)$.\end{remark}

%\begin{remark} (Controlling $J_j(K^\star) - J_j(K^\star_j)$) We emphasize that similarly to \citep[Theorem 4]{toso2024meta} we can also control the gap between the optimal meta-controller $K^\star$ and any task-specific optimal controller $K^\star_j$. In particular, $J_j(K^\star) - J_j(K^\star_j) =\mathcal{O}( f^2(\epsilon_{\text{het}}))$, which guarantees that $K^\star$ is close to $K^\star_j$ up to a task heterogeneity bias. \end{remark}

\subsection{Proof of Corollary \ref{corollary: sample complexity global} (Sample Complexity)} \label{appendix: sample complexity global}

We let $\bar{\mathcal{S}}^{\mathcal{S}}_c$ and $\bar{\mathcal{S}}^{\mathcal{M}}_c$ to denote the total number of samples required in Algorithm \ref{alg:metatrain}, for the MAML-LQR problem, to learn a meta-controller $K_N$ that is $\epsilon$-close to any task-specific optimal controller up to a heterogeneity bias, with and without task selection, respectively. In particular, $\bar{\mathcal{S}}^{\mathcal{S}}_c = \mathcal{O}\left(LNn_s\right) + \mathcal{O}\left(Mn_s\right)$ and $\bar{\mathcal{S}}^{\mathcal{M}}_c = \mathcal{O}\left(MNn_s\right)$. In addition, to guarantee \eqref{eq:global convergence analysis} with high probability, we need $n_s = \mathcal{O}(d^2/\epsilon^2)$ samples in the zeroth-order gradient estimation. Therefore, $\bar{\mathcal{S}}^{\mathcal{M}}_c = \mathcal{O}\left(MNn_s\right) = \mathcal{O}((d^2/\epsilon^{2})\log(1/\epsilon))$, whereas
$\bar{\mathcal{S}}^{\mathcal{S}}_c = \frac{L}{M}\mathcal{O}\left((d^2/\epsilon^{2})\log(1/\epsilon)\right) + \mathcal{O}\left(d^2/\epsilon^{2}\right),$
then, for a sufficiently large amount of tasks in the task pool $\mathcal{M}$ (e.g., scaling as $M = \mathcal{O}(\log(1/\epsilon))$) and a sufficiently small amount of tasks in the subset $\mathcal{S}$, i.e., $L = \mathcal{O}(1)$, the task selection may benefit from a sample complexity reduction of a factor of up to $\textcolor{blue}{\mathcal{O}(\log(1/\epsilon))}$ when compared to the setting without task selection.

\subsection{Stability Analysis} \label{appendix: stability analysis}

We now proceed to demonstrate that $K_n \in \mathcal{G}$, for any $n = \{0,1,\ldots,N-1\}$ of Algorithm \ref{alg:metatrain}. Let us first recall that by Assumption \ref{assumption:initial_stabilizing_K0}, the initial meta-controller is stabilizing, i.e., $K_0 \in \mathcal{G}$. Then, we can first show that $\bar{K}_0 = K_0 - \eta_{\text{inn}}\widehat{\nabla} J_j(K_0)$ and $K_1 = K_0 - \eta_{\text{out}}\nabla_{\mathcal{S}}J(K_0)$ never leaves $\mathcal{G}$, with high probability, if the number of samples $n_s$, smoothing radius $r$, inner and outer step-sizes $\eta_{\text{inn}}$, $\eta_{\text{out}}$, and heterogeneity $f(\epsilon_{\text{text}})$ are set accordingly. Finally, we can use an induction step to extend the same for any iteration $n$. By the gradient smoothness in Lemma \ref{lemma:LQR} we can write
\begin{align*}
    J_j(\bar{K}_0) - J_j(K_0) &\leq  -\frac{\eta_{\text{inn}}}{4}\|{\nabla} J_j(K_0)\|^2_F + \frac{3\eta_{\text{inn}}}{4}\|\widehat{\nabla} J_j(K_0) - \nabla J_j(K_0)\|^2_F,
\end{align*}
where $\eta_{\text{out}} \leq \frac{1}{4\psi}$. Then, by using the gradient dominance property we have that

\begin{align*}
    J_j(\bar{K}_0) - J_j(K^\star_j) &\leq  \left( 1 -\frac{\lambda_j\eta_{\text{inn}}}{4}\right)\Delta^{(j)}_0 + \frac{3\eta_{\text{inn}}}{4}\|\widehat{\nabla} J_j(K_0) - \nabla J_j(K_0)\|^2_F,
\end{align*}
where $\|\widehat{\nabla} J_j(K_0) - \nabla J_j(K_0)\|^2_F$ corresponds to the zeroth-order gradient estimation error at $K_0$. As well-established in \citep{toso2024meta,toso2024asynchronous} and also discussed previously in this work, the zeroth-order estimation error can be made arbitrarily small, for instance, $\|\widehat{\nabla} J_j(K_0) - \nabla J_j(K_0)\|_F \leq \sqrt{\frac{\lambda_j\Delta^{(j)}_0}{6}}$, for $n_s = \mathcal{O}\left(\frac{6\psi}{\lambda_j}\right)$ and $r = \mathcal{O}\left(\sqrt{\frac{\lambda_j\Delta^{(j)}_0}{6}}\right)$. Then, for all tasks $T_j \sim \mathcal{P}(\mathcal{T})$ we have
\begin{align*}
    J_j(\bar{K}_0) - J_j(K^\star_j) &\leq  \left( 1 -\frac{\lambda_j\eta_{\text{inn}}}{8}\right)\Delta^{(j)}_0,
\end{align*}
which implies that $K_0 \in \mathcal{G}$ (i.e., see Definition \ref{def:stabilizing_set}). We proceed to show that $K_1 \in \mathcal{G}$. To do so, we use again the gradient smoothness in Lemma \ref{lemma:LQR} to write

\begin{align*}
    J_j(K_1) &- J_j(K_0) \leq  -\frac{\eta_{\text{out}}}{4}\|{\nabla} J_j(K_0)\|^2_F + \frac{3\eta_{\text{out}}}{4}\|\nabla_{\mathcal{S}}J(K_0) - \nabla J_j(K_0)\|^2_F\\
    &\stackrel{(i)}{\leq} -\frac{\eta_{\text{out}}}{4}\|{\nabla} J_j(K_0)\|^2_F + \frac{3\eta_{\text{out}}}{2}\|\nabla_{\mathcal{S}}J(K_0) - \nabla J(K_0)\|^2_F + 3\eta_{\text{out}}f^2(\epsilon_{\text{het}}) + 3\eta_{\text{out}}\eta^2_{\text{inn}}\phi^2,
\end{align*}
where $(i)$ is due to Young's inequality \eqref{eq:youngs} and $$\|\nabla J(K_0) - J_j(K_0)\|^2_F \leq \mc{E}_i\|\nabla J_i(K_0 - \eta_{\text{inn}}\nabla J_i(K_0)) - \nabla J_j(K_0)\|^2_F \leq 2\eta^2_{\text{inn}}\phi^2 + 2f^2(\epsilon_{\text{het}}).$$

Then, by the gradient dominance property we have that
\begin{align*}
    J_j(K_1) - J_j(K^\star_j) &\leq \left(1 -\frac{\lambda_j\eta_{\text{out}}}{4}\right) \Delta^{(j)}_0 + \frac{3\eta_{\text{out}}}{2}\|\nabla_{\mathcal{S}}J(K_0) - \nabla J(K_0)\|^2_F\\
    &+ 3\eta_{\text{out}}f^2(\epsilon_{\text{het}}) + 3\eta_{\text{out}}\eta^2_{\text{inn}}\phi^2,
\end{align*}
where the gradient estimation error $\|\nabla_{\mathcal{S}}J(K_0) - \nabla J(K_0)\|_F$ satisfy

\begin{align*}
    \|\nabla_{\mathcal{S}}J(K_0) - \nabla J(K_0)\|_F \leq \zeta +\frac{1}{M}\sum_{j \in \mathcal{M}} \min_{i \in \mathcal{S^\star}}\xi_{i,j} \leq \zeta + f(\epsilon_\text{het}),
\end{align*}
with probability $1-\delta$, for $n_s$ and $r$ satisfying \eqref{eq: ns and r} with $\zeta$ in lieu of $\epsilon$. Then, we can write 
\begin{align*}
    J_j(K_1) - J_j(K^\star_j) &\leq \left(1 -\frac{\lambda_j\eta_{\text{out}}}{4}\right) \Delta^{(j)}_0 + 3\eta_{\text{out}}\zeta^2 + 6\eta_{\text{out}}f^2(\epsilon_{\text{het}}) + 3\eta_{\text{out}}\eta^2_{\text{inn}}\phi^2\\
    &\stackrel{(i)}{\leq} \left(1 -\frac{\lambda_j\eta_{\text{out}}}{8}\right) \Delta^{(j)}_0
\end{align*}
where $(i)$ follows from $\zeta = \sqrt{\frac{\lambda_j\Delta^{(j)}_0}{72}}$, $f(\epsilon_{\text{het}}) \leq \frac{\lambda_j \Delta^{(j)}_0}{72}$, and $\eta_{\text{inn}} \leq \sqrt{\frac{\lambda_j\Delta^{(j)}_0}{72 \phi^2}}$, which implies that $K_1 \in \mathcal{G}$. Therefore, we define our base case and inductive hypothesis as follows:

\begin{align*}
    \textbf{Base case:} \;\ J_j(\bar{K}_0) - J_j(K^\star_j) &\leq  J_j(K_0) - J_j(K^\star_j),\\
    J_j(K_1) - J_j(K^\star_j) &\leq  J_j(K_0) - J_j(K^\star_j)
\end{align*}

\begin{align*}
    \textbf{Inductive hypothesis:} \;\ J_j(\bar{K}_n) - J_j(K^\star_j) &\leq  J_j(K_0) - J_j(K^\star_j),\\
    J_j(K_n) - J_j(K^\star_j) &\leq  J_j(K_0) - J_j(K^\star_j)
\end{align*}
which can be used along with the aforementioned conditions on the number of samples, smoothing radius, step-sizes and heterogeneity, to write 
\begin{align*}
    J_j(K_{n+1}) - J_j(K^\star_j) &\leq \left(1 -\frac{\lambda_j\eta_{\text{out}}}{4}\right) \Delta^{(j)}_n + 3\eta_{\text{out}}\zeta^2 + 6\eta_{\text{out}}f^2(\epsilon_{\text{het}}) + 3\eta_{\text{out}}\eta^2_{\text{inn}}\phi^2\\
    &\stackrel{(i)}{\leq} \left(1 -\frac{\lambda_j\eta_{\text{out}}}{8}\right) \Delta^{(j)}_0 \leq \Delta^{(j)}_0,
\end{align*}
which guarantees that Algorithm \ref{alg:metatrain} produces MAML stabilizing controllers with high probability.

\end{document}